\documentclass[11pt]{article}
\input{epsf.sty}
\usepackage{hyperref}
\usepackage[final]{epsfig}
\usepackage{color}
\usepackage{dsfont}

\usepackage{amsmath, amsthm, amssymb}
\newtheorem {thm}{Theorem}[section]
\newtheorem {prop}[thm]{Proposition}
\newtheorem {lem}[thm]{Lemma}
\newtheorem {cor}[thm]{Corollary}
\newtheorem {defn}[thm]{Definition}

\newtheorem {rem}[thm]{Remark}

\def\Cox{\hfill \Box}

\def\N{{\Bbb N}}
\def\TT{{\Bbb T}}
\def\Z{{\Bbb Z}}
\def\R{{\Bbb R}}
\def\P{{\Bbb P}}

\def\E{{\Bbb E}}

\newcommand{\bzd}{\ensuremath{{(\Z^d)^*}}}
\newcommand{\bzdl}{\ensuremath{{\Lambda^*}}}

\newcommand{\C}{{\mathcal C}}
\newcommand{\A}{{\mathcal A}}
\newcommand{\ormd}{\mathrm{d}}
\newcommand{\rmd}{\,\mathrm{d}}

\newcommand{\shift}{\mbox{\scriptsize shift}}

\def\0{{\bf 0}}

\def\var{{\rm var}}

\newcommand{\zd}{{{\mathbb Z}^d}}
\newcommand{\RR}{\ensuremath{\mathbb{R}}}

\def\phi{\varphi}

\def\g{\gamma}

\def\k{\kappa}

\def\x{\xi}

\def\o{\omega}

\def\L{\Lambda}

\def\T{\T}

\def\AA{{\cal A}}
\def\C{{\cal C}}
\def\GG{{\cal G}}
\def\PP{{\cal P}}
\def\NN{{\cal N}}
\def\FF{{\cal F}}
\def\HH{{\cal H}}
\def\TT{{\cal T}}
\def\Norm{\text{Norm}\,}
\def\Const{\text{Const}\,}
\def\const{\text{const}\,}

\newcommand{\cov}{{\mbox{\rm\normalfont cov}}}
\newcommand{\Cov}{\ensuremath{\Bbb C\mbox{\rm\normalfont ov}}}
\newcommand{\Var}{\ensuremath{\Bbb V\mbox{\rm\normalfont ar}}}

\begin{document}

\title{Uniqueness of gradient Gibbs measures with disorder
} 

\author{
Codina Cotar
\footnote{
University College London, Statistical Science Department,
London,
United Kingdom,
\texttt{c.cotar@ucl.ac.uk}}
\, and
Christof K\"ulske
\footnote{
Ruhr-University of Bochum, Fakult\"at f\"ur Mathematik, Bochum,
Germany,
\texttt{Christof.Kuelske@ruhr-uni-bochum.de}}
}

\newcommand{\CC}[1]{{\color{blue} #1}}
\newcommand{\DE}[1]{{\color{red} #1}}
\newcommand{\CK}[1]{{\color{green} #1}}

\def\AA{{\cal A}}
\def\C{{\cal C}}
\def\DD{{\cal D}}
\def\GG{{\cal G}}
\def\PP{{\cal P}}
\def\NN{{\cal N}}
\def\FF{{\cal F}}
\def\HH{{\cal H}}
\def\TT{{\cal T}}
\def\Norm{\text{Norm}\,}
\def\Const{\text{Const}\,}
\def\const{\text{const}\,}

\maketitle

\begin{abstract}
We consider - in a {uniformly strictly convex potential regime} - two versions of random gradient models with disorder. 
In model (A) the interface feels a bulk term of random fields 
while in model (B) the disorder enters though the potential acting on the gradients. We assume a \textit{general distribution} on the disorder with \textit{uniformly-bounded} finite second moments. 

It is well known that for gradient models without disorder
there are no Gibbs measures in infinite-volume in dimension
$d = 2$, while there are shift-invariant gradient Gibbs measures describing an
infinite-volume distribution for the gradients of the field, as was
shown by Funaki and Spohn in \cite{FS}. Van Enter and K\"ulske proved in \cite{EK} that adding a disorder term as in model (A) prohibits
the existence of such gradient Gibbs measures for general interaction
potentials in $d = 2$. In \cite{CK} we proved the existence of shift-covariant random gradient Gibbs measures for model (A) when $d\geq 3$, the disorder is i.i.d  and has mean zero, and for model (B) when $d\geq 1$ and the disorder has a stationary distribution. 

In the present paper, we prove \textit{existence and uniqueness} of shift-covariant random gradient Gibbs measures with a given \textit{expected tilt} $u\in\R^d$ and with the corresponding annealed measure being \textit{ergodic}: for model (A) when $d\geq 3$ and the disordered random fields are \textit{i.i.d.} and \textit{symmetrically-distributed}, and for model (B) when $d\geq 1$ and for \textit{any stationary} disorder-dependence structure. We also compute for both models for any gradient Gibbs measure constructed as in \cite{CK}, when the disorder is \textit{i.i.d.} and its distribution satisfies a \textit{Poincar\'e inequality} assumption, the \textit{optimal decay of covariances} with respect to the \textit{averaged-over-the-disorder} gradient Gibbs measure.
\end{abstract}

\bigskip 


\section{Introduction} \label{sect:intro}
Phase separation in $\RR^{d+1}$ can be described by effective interface models for the study of phase boundaries at a mesoscopic level
in statistical mechanics. Interfaces are sharp boundaries which separate the different regions of space occupied by different phases. In this class of models, the interface is modeled as the graph of a random
function from $\Z^d$ to $\Z$ or to $\RR$ (discrete or continuous effective interface models). For background and earlier results on continuous and discrete interface models without disorder see for example \cite{BY}, \cite{CD}, \cite{CDM}, \cite{DGI}, \cite{FL}, \cite{FS}, \cite{gos} and references therein. In our setting, we will consider the case of continuous interfaces with disorder as introduced and studied previously in \cite{EK} and \cite{KO}. Note also that discrete interface models in the presence of
disorder have been studied for example in \cite{BK1} and \cite{BK2}. 

There is some similarity between models of continuous interfaces and models of rotators ($S^1$-valued spins) which interact 
via a spin-rotation invariant ferromagnetic interaction.  
It is a classical result of mathematical physics that, at low enough temperatures, 
there is a continuous symmetry breaking and ferromagnetic order in these rotator models for space dimensions $d\geq 3$, at (for Lebesgue) a.e. temperature, 
see \cite{FSS} and \cite{P}. Generally speaking, adding disorder to a model tends to destroy the non-uniqueness of Gibbs measures, 
and to destroy \textit{order}, for the precise statements see \cite{AW}. Indeed the non-existence results for interfacial states of  \cite{BK1} and \cite{EK}  rely on suitable adaptations of this method.   

Nevertheless, there are striking examples where disorder acts in an opposite way: 
Non-uniqueness of the Gibbs measure and a new type of ordering 
can even be created by the introduction of quenched randomness of a random field type. Such an order-by-disorder 
mechanism was proved to happen in the rotator model in the presence of  a uni-axial random field, see \cite{C1} and \cite{C2}. In this model the rotators tend to align in a plane perpendicular to the axis of the external fields. 
Heuristically it seems that the mechanism for such a random-field-induced order should remain particular to models of rotators, since 
the interplay of disorder, interaction, and boundedness of spins is crucial.  

However, this example underlines the subtlety of the uniqueness issue for continuous models which are subjected to random fields in general.

\subsection{Our models}

We will introduce next our two models of interest. 

In our setting, the fields $\phi(x)\in\RR$ represent height variables
of a random interface at the sites $x\in\zd$. Let $\Lambda$ be a finite set in $\Z^d$ with boundary 
\begin{eqnarray}
\partial\Lambda:=\{x\notin\Lambda,~||x-y||=1~\mbox{for some}~y\in\Lambda\},~\mbox{where}~\|x-y\|_=\sum^d_{i=1}|x_i-y_i|.
\end{eqnarray}
On the boundary we set a boundary condition $\psi$ such that $\phi(x)=\psi(x)$ for $x\in\partial\Lambda$. Let $(\Omega, {\cal F},\P)$ be a probability 
space; this is the probability space of the disorder, which will be introduced below. We denote by the symbol $\E$ the expectation w.r.t $\P$, by $\Var$ the variance w.r.t. $\P$ and by $\Cov$ the covariance w.r.t $\P$.

Our two models are given in terms of the \textit{finite-volume Hamiltonian} on $\Lambda$.
\begin{enumerate}
\item [(A)] For model A the Hamiltonian is
\begin{eqnarray}
\label{eqn00}
H_{\Lambda}^\psi[\xi](\phi):=\frac{1}{2}\sum_{x, y\in\Lambda\atop |x-y|=1}V(\phi(x)-\phi(y))+\sum_{x\in\Lambda,y\in\partial\Lambda\atop |x-y|=1}V(\phi(x)-\psi(y))+\sum_{x\in\Lambda}\xi(x)\phi(x),
\end{eqnarray}
where the random fields $(\xi(x))_{x\in\Z^d}$ are assumed to be \textit{i.i.d.} real-valued random variables, with \textit{finite non-zero second moments}.
The disorder configuration $(\xi(x))_{x\in\Z^d}$ denotes an
arbitrary fixed configuration of external fields, modeling a ``quenched'' (or
frozen) random environment. We assume that $V\in C^2(\Bbb R)$ is an even function 
such that there exist $0<C_1<C_2$ with
\begin{equation}
\label{tag22}
C_1\le V''(s)\le C_2 ~\mbox{for all}~ s\in\RR.
\end{equation} 
\item [(B)] For each bond $(x,y)\in\Z^d\times\Z^d,|x -y|=1$, we define the measurable map $V_{(x,y)}^{\omega}(s):(\omega,s)\in\Omega\times\RR\rightarrow\RR$. Then $V_{(x,y)}^\omega$ is a random real-valued function. 
Assume that $V_{(x,y)}^\omega\in C^2(\R)$ have \textit{uniformly-bounded finite second moments} and jointly \textit{stationary} distribution.
We also assume that for some given $0<C_{1, (x,y)}^\omega<C_{2,(x,y)}^\omega,\omega\in\Omega$, with $0<\inf_{(x,y)}\E\big(C_{1, (x,y)}^\omega\big)<\sup_{(x,y)}\E\big(C_{2, (x,y)}^\omega\big)<\infty$, $V_{(x,y)}^\omega$ obey for $\P$-almost every $\omega\in\Omega$ the following bounds, uniformly in the bonds $(x,y)$
\begin{equation}
\label{tag23}
{C}^\omega_{1,(x,y)}\le (V_{(x,y)}^\omega)''(s)\leq {C}^\omega_{2,(x,y)} ~\mbox{for all}~s\in\RR.
\end{equation}
We set the further condition that for each fixed $\omega\in\Omega$ and for each bond $(x,y)$,  $V_{(x,y)}^\omega\in C^2(\Bbb R)$ is an even function. Then for model B we define the Hamiltonian for each fixed $\omega\in\Omega$ by
\begin{eqnarray}
\label{eqn000}
H_{\Lambda}^\psi[\omega](\phi):=\frac{1}{2}\sum_{x,y\in\Lambda,|x-y|=1}V_{(x,y)}^\omega(\phi(x)-\phi(y))+\sum_{x\in\Lambda,y\in\partial\Lambda,|x-y|=1}V_{(x,y)}^\omega(\phi(x)-\psi(y)).
\end{eqnarray}
\end{enumerate}
For our second main result for both models A and B, we will work under the following slightly more restrictive Poincar\'e inequality assumption on the distribution $\gamma$ of the disorder $\xi(0)$, (respectively of $V^\omega_{(0,e_1)}$):
There exists $\lambda>0$ such that for all smooth enough real-valued functions $f$ on $\Omega$, we have for the probability measure $\gamma$
\begin{equation}
\label{gap}
\lambda\var_{\gamma}(f)\le\int |\nabla f|^2\rmd\gamma,
\end{equation}
where $|\nabla f|$ is the Euclidean norm of the gradient of f and $\var_\gamma$ is the variance with respect to $\gamma$. By smooth, we
understand in the above enough regularity in order that the various
expressions we are dealing with are well defined and finite. 
Known examples where the Poincar\'e inequality holds have been described by the so-called Bakry-Emery criterion \cite{BAE}, which involves log-concavity conditions on the measure rather than on its density. For further explicit assumptions on $\gamma$ such that (\ref{gap}) holds, see for instance \cite{Led} or (for a large class of non-convex potentials) Theorem 3.8 from \cite{Mil}.

\begin{rem}  Our model B with uniformly strictly convex potentials is the gradient model analogue of the random conductance model with uniform ellipticity condition. See, for example, \cite{BB} for an extensive review on the random conductance model and its connection to the gradient model.\end{rem}


The two models above
are prototypical ways to add randomness which preserves the gradient structure, i.e., the Hamiltonian depends only on the gradient field $(\phi(x)-\phi(y))_{x,y\in\Z^d, |x-y|=1}$. Note that for $d=1$ 
our interfaces can be used to model a polymer chain, see for example \cite{denholl}. Disorder in the Hamiltonians models impurities in the physical system. Models A and B can be regarded as modeling two different types of impurities, one affecting the interface height, the other affecting the interface gradient. 

The rest of the introduction is structured as follows: in Subsection \ref{AG} we define in detail the notions of finite-volume and infinite-volume (gradient) Gibbs measures for model A, in Subsection \ref{BG} we sketch the corresponding notions for 
model B, and in Subsection \ref{MR} we present our main results and their connection to the existing literature.

\subsection{Gibbs measures and gradient Gibbs measures for model A}
\label{AG}
\subsubsection{$\phi$-Gibbs measures}

Let $C_b(\R^{\zd})$ denote the set of continuous and bounded functions on $\R^{\zd}$. The functions considered are functions of the interface configuration $\phi$, and
continuity is with respect to each coordinate $\phi(x),x\in\Z^d,$ of the interface. For a finite region $\Lambda\subset\Z^d$, let $\rmd\phi_{\Lambda}:=\prod_{x\in\Lambda}\rmd\phi(x)$ be the Lebesgue measure over $\R^{\Lambda}$.

\noindent Let us
first consider model A only, and let us define the $\phi$-Gibbs measures for {\bf
fixed} disorder $\xi$.
\begin{defn} {\bf(Finite-volume $\phi$-Gibbs measure)}
\label{gibbs0}
For a finite region $\Lambda\subset\Z^d$, {\em the finite-volume Gibbs measure $\nu_{\Lambda,\psi}[\xi]$ on $\R^{\Z^d}$} with given Hamiltonian $H[\xi]:=(H_{\Lambda}^\psi[\xi])_{\Lambda\subset\zd, \psi \in
\RR^{\zd}}$, with boundary condition $\psi$ for the field of height variables $(\phi(x))_{x\in\Z^d}$ over $\Lambda$, and with a fixed disorder configuration $\xi$, is defined by
\begin{equation}
\label{tag0'}
\nu_{\Lambda}^\psi[\xi ](\ormd\phi):=\frac{1}{Z_{\Lambda}^\psi[\xi]}\exp\left\{-H_{\Lambda}^\psi[\xi](\phi)\right\}\rmd\phi_\Lambda\delta_\psi(\ormd\phi_{{\Z}^d\setminus\Lambda}). 
\end{equation}
where
$$Z_{\Lambda}^\psi[\xi]:=\int_{{\R}^{\Z^d}}\exp\left\{-H_{\Lambda}^\psi[\xi](\phi)\right\}\rmd\phi_\Lambda\delta_\psi(\ormd\phi_{\Z^d\setminus\Lambda})$$
and
$$\delta_\psi(\ormd\phi_{\Z^d\setminus\Lambda}):=\prod_{x\in\Z^d\setminus\Lambda}\delta_{\psi(x)}(\ormd\phi(x)).$$ 
\end{defn}
It is easy to see that the conditions on $V$ guarantee the finiteness of the integrals appearing in (\ref{tag0'}) for
all arbitrarily fixed choices of $\xi$.
\begin{defn}  {\bf ($\phi$-Gibbs measure on $\zd$)}
\label{gibbs}
The probability measure $\nu[\xi]$ on $\RR^{\zd}$ is called an {\em (infinite-volume) Gibbs measure} for the $\phi$-field with given Hamiltonian $H[\xi]:=(H_{\Lambda}^\psi[\xi])_{\Lambda\subset\zd, \psi \in
\RR^{\zd}}$ ($\phi$-Gibbs measure for short), if it satisfies the DLR equation
\begin{equation}
\label{dlrgibbs}
\int\nu[\xi](\ormd\psi)\int \nu_{\Lambda}^\psi[\xi](\ormd\phi)F(\phi)=\int\nu[\xi](\ormd\phi)F(\phi),
\end{equation}
for every finite $\Lambda\subset\zd$ and for all $F\in C_b(\R^{\zd})$.
\end{defn}
We discuss next the case of interface models without disorder, that is, with $\xi(x)=0$ for all $x\in \Z^d$ in model A. 
Let $\nu^\psi_\Lambda[\xi=0],\L\in\Z^d$, denote the finite-volume Gibbs measure for $\L$ and with boundary condition $\psi$. Then an infinite-volume Gibbs measure $\nu[\xi=0]$ exists under the conditions $V(s)\ge As^2+B$ and $V''(s)\le C_2$, $A, C_2>0, B\in\R, s\in\R,$ 
only when $d\ge3$, but not for $d=1,2$, where the field "delocalizes" as $\Lambda\nearrow\Z^d$ (see \cite{FP}). 

In the case of interfaces with disorder as in model A, it has been proved in \cite{KO} that the $\phi$-Gibbs measures do not exist when $d=2$. A similar argument as in \cite{KO} can be used to show that $\phi$-Gibbs measures do not exist for model A when $d=1$.

\subsubsection{$\nabla\phi-$Gibbs Measures}
\label{nablagibbs}
We note that the Hamiltonian $H_\L^\psi[\xi]$ in model A, respectively $H_\L^\psi[\omega]$ in model B, changes only by a configuration-independent
constant under the joint shift $\phi(x)\rightarrow\phi(x)+c$ of all height variables $\phi(x),x\in\Z^d,$
with the same $c\in\RR$. This holds true for any fixed configuration $\xi$, respectively $\omega$. Hence, finite-volume
Gibbs measures transform under a shift of the boundary condition
by a shift of the integration variables. Using this invariance under height
shifts we can lift the finite-volume measures to measures on gradient configurations, i.e., configurations of height differences across
bonds,
defining the gradient finite-volume Gibbs measures. Gradient Gibbs measures have the advantage that they may exist, even
in situations where the Gibbs measure does not. Note that the concept of $\nabla\phi-$ measures is general and does not refer only to the disordered models. For example, in the case of interfaces without disorder $\nabla\phi$-Gibbs measures exist for all $d\ge 1$. 

We next introduce the \textit{bond variables on $\zd$}. Let 
$$\bzd:=\{b=(x_b,y_b)~|~x_b,y_b\in\zd,\|x_b-y_b\|=1,b~\mbox{directed from}~x_b~\mbox{to}~y_b\},$$ 
where $\|x\|=\max_{1\le i\le d} |x_i|,$ for $x=(x_1,\ldots,x_d)\in\zd$; note that each undirected bond appears twice in $\bzd$. 	We define
$$\bzdl:=\bzd\cap (\Lambda\times\Lambda)~~\mbox{and}~~\partial\bzdl:=\{b=(x_b,y_b)~|~x_b\in\zd\setminus\Lambda,y_b\in\Lambda,\|x_b-y_b\|=1\}.$$
For $\phi=(\phi(x))_{x\in\zd}$ and $b=(x_b,y_b)\in\bzd$, we define the \textit{height differences} $\nabla\phi(b):=\phi(y_b)-\phi(x_b)$.
The height variables $\phi=\{\phi(x):x\in\zd\}$ on $\zd$ automatically determine a field of height differences
$\nabla\phi=\{\nabla\phi(b):b\in\bzd\}$. One can therefore consider the distribution $\mu$ of $\nabla\phi$-fields
under the $\phi$-Gibbs measure $\nu$. We shall call $\mu$ the $\nabla\phi$-Gibbs measure. In fact, it is possible to define
the $\nabla\phi$-Gibbs measures directly by means of the DLR equations and, in this sense, $\nabla\phi$-Gibbs measures exist for
all dimensions $d\ge1$.

A sequence of bonds $\C=\{b^{(1)},b^{(2)},\ldots,b^{(n)}\}$ is called a \textit{chain} connecting $x$ and $y$, $x,y\in\zd$, if $x_{b_1}=x,y_{b^{(i)}}=x_{b^{(i+1)}}$ for $1\le i\le n-1$ and $y_{b^{(n)}}=y$. The chain is called a \textit{closed loop} if $y_{b^{(n)}}=x_{b^{(1)}}$. A \textit{plaquette} is a closed loop $\A=\{b^{(1)},b^{(2)},b^{(3)},b^{(4)}\}$ such that $\{x_{b^{(i)}},i=1,\ldots,4\}$ consists of 
$4$ different points. 

The field $\eta=\{\eta(b)\}\in\RR^{\bzd}$ is said to satisfy \textit{the plaquette condition} if
\begin{eqnarray}
\eta(b)=-\eta(-b)~\mbox{for all}~b\in\bzd~\mbox{and}~\sum_{b\in\A}\eta(b)=0~\mbox{for all plaquettes}~\A~\mbox{in}~\zd,
\end{eqnarray}
where $-b$ denotes the reversed bond of $b$. Let 
\begin{equation}
\chi=\{\eta\in\RR^{(\zd)^*}~\mbox{which satisfy the plaquette 
condition}\}
\end{equation}
and let $L_r^2, r>0$, be the set of all $\eta\in\RR^{\bzd}$ such that
$$|\eta|^2_r:=\sum_{b\in\bzd}|\eta(b)|^2e^{-2r\|x_b\|}<\infty.$$
We denote $\chi_r=\chi\cap L_r^2$ equipped with the norm $|\cdot|_r$. For $\phi=(\phi(x))_{x\in\zd}$ and $b\in\bzd$, we define $\eta(b):=\nabla\phi(b)$. Then $\nabla\phi=\{\nabla\phi(b):b\in\bzd\}$ satisfies the plaquette condition. Conversely, the heights $\phi^{\eta,\phi(0)}\in\RR^\zd$ can be constructed from height differences $\eta$ and the
height variable $\phi(0)$ at $x=0$ as 
\begin{equation}
\label{19}
\phi^{\eta,\phi(0)}(x):=\sum_{b\in\C_{0,x}}\eta(b)+\phi(0),
\end{equation} 
where $\C_{0,x}$ is an arbitrary chain connecting $0$ and $x$. Note that $\phi^{\eta,\phi(0)}$ is well-defined if
$\eta=\{\eta(b)\}\in\chi$.

Let $C_b(\chi)$ be the set of continuous and bounded functions on $\chi$, where the continuity is with respect to each bond variable $\eta(b),b\in\bzd$.
\begin{defn}{\bf(Finite-volume $\nabla\phi$-Gibbs measure)}
\label{finvolgrad}
The {\em finite-volume $\nabla\phi$-Gibbs measure} in $\Lambda$ (or more precisely, in $\bzdl$) with given Hamiltonian $H[\xi]:=(H_{\Lambda}^ \rho[\xi])_{\Lambda\subset\zd,\,\rho \in
\chi}$, with boundary condition $\rho\in\chi$ and with fixed disorder configuration $\xi$, is a probability measure $\mu_{\Lambda}^\rho[\xi]$ on $\chi$ such that for all $F\in C_b(\chi)$, we have
\begin{equation}
\label{fingradgibbs}
\int_{\chi}\mu_{\Lambda}^\rho[\xi](\ormd\eta)F(\eta)=\int_{\R^{\zd}}\nu_{\Lambda}^\psi[\xi](\ormd\phi)F(\nabla\phi),
\end{equation}
where $\psi$ is any field configuration whose gradient field is $\rho$. 
\end{defn}
\noindent We are now ready to define the main object of interest of this paper: the {\bf random} (gradient) Gibbs measures. 
\begin{defn} {\bf ($\nabla\phi$-Gibbs measure on $(\zd)^*$)} 
\label{nablaphigib}
The probability measure $\mu[\xi]$ on $\chi$ is called an {\em (infinite-volume) gradient Gibbs measure} with given Hamiltonian $H[\xi]:=(H_{\Lambda}^\rho[\xi])_{\Lambda\subset\zd, \rho \in\chi}$ ($\nabla\phi$-Gibbs measure for short),
if it satisfies the DLR equation
\begin{equation}
\label{dlrgrad}
\int\mu[\xi](\ormd\rho)\int \mu_{\Lambda}^\rho[\xi](\ormd\eta)F(\eta)=\int\mu[\xi](\ormd\eta)F(\eta),
\end{equation}
for every finite $\Lambda\subset\zd$ and for all $F\in C_b(\chi)$.
\end{defn}
\begin{rem}
Throughout the rest of the paper, we will use the notation $\phi,\psi$ to denote height variables and $\eta,\rho$ to denote gradient variables.
\end{rem}
For $v\in\zd$, we define the shift operators: $\tau_{v}$ for the heights by $(\tau_{v}\phi)(y):=\phi(y-v)~\mbox{for}~y\in\zd~\mbox{and}~\phi\in\RR^{\zd}$, $\tau_{v}$ for the bonds by $(\tau_{v}\eta)(b):=\eta(b-v)$ for $b\in\bzd~\mbox{and}~\eta\in\chi$, and $\tau_v$ for the disorder configuration by $(\tau_v\xi)(y):=\xi(y-v)$ for $y\in\zd$ and $\xi\in\RR^{\zd}$.
\begin{defn} {\bf(Translation-covariant random (gradient) Gibbs measures for model A)}
\label{shiftcov1}
A measurable map $\xi\rightarrow\nu[\xi]$ is called {\em a translation-covariant random Gibbs
measure} if $\nu[\xi]$ is a $\phi$-Gibbs measure for $\P$-almost every $\xi$, and if 
$$\int\nu[\tau_v\xi](\ormd\phi)F(\phi)=\int\nu[\xi](\ormd\phi)F(\tau_v\phi),$$
for all $v\in\zd$ and for all $F\in C_b(\R^{\zd})$.

To define the notion of measurability for a measure-valued function we use the 
evaluation sigma-algebra in the image space, which is the smallest sigma-algebra such that
the evaluation maps $\mu\mapsto \mu(A)$ are measurable for all events $A$ (for details, see page 129 from Section 7.3 on the extreme decomposition in \cite{giorgii}).


A measurable map $\xi\rightarrow\mu[\xi]$ is called {\em a translation-covariant random gradient Gibbs measure} if $\mu[\xi]$ is a $\nabla\phi-$ Gibbs measure for $\P$-almost every $\xi$, and if 
$$\int\mu[\tau_v\xi](\ormd\eta)F(\eta)=\int\mu[\xi](\ormd\eta)F(\tau_v\eta),$$
for all $v\in\zd$ and for all $F\in C_b(\chi)$.
\end{defn}
The above notion generalizes the
notion of a translation-invariant (gradient) Gibbs measure to the set-up of disordered
systems.
\begin{rem}
Throughout the paper, we will use the notation $\nu_\Lambda$, respectively $\nu$, to denote a finite-volume, respectively the corresponding infinite-volume, Gibbs measure, and the notation $\mu_\Lambda$, respectively $\mu$, to denote a finite-volume, respectively the corresponding infinite-volume, gradient Gibbs measure.
\end{rem}

\subsection{Gibbs measures and gradient Gibbs measures for model B}
\label{BG}

The notions of finite-volume (gradient) Gibbs measure and infinite-volume (gradient) Gibbs measure for model B can be defined similarly as for model A, with  $(V^\omega_{(x,y)})_{(x,y)\in\zd\times\zd},\omega\in\Omega$, playing a similar role to $\xi\in\RR^\zd$, and with $\omega$ replacing $\xi$ in Definitions \ref{gibbs0}-\ref{nablaphigib}. Once we specify the action of the shift
map $\tau_v$ in this case, we can also define the notion of translation-covariant random (gradient) Gibbs measure, with $\omega\in\Omega$ replacing $\xi\in\RR^\zd$ in Definition \ref{shiftcov1}.

Let $\tau_v,v\in\Z^d,$ be a shift-operator and let $\omega\in\Omega$ be fixed. We 
will denote by $\nu[\tau_v\omega]$ the infinite-volume Gibbs measure with given Hamiltonian $\bar{H}[\omega](\phi):=\left(H_\L^\psi[\omega](\tau_v\phi)\right)_{\L\subset\Z^d,\psi\in\RR^{\Z^d}}$. This means that we shift the field of disorded potentials on bonds from $V_{(x,y)}^\omega$ to $V_{(x+v,y+v)}^\omega$. Similarly, we 
will denote by $\mu[\tau_v\omega]$ the infinite-volume gradient Gibbs measure with given Hamiltonian $\bar{H}[\omega](\eta):=\left(H_\L^\rho[\omega](\tau_v\eta)\right)_{\L\subset\Z^d,\rho\in\RR^{\bzd}}$. 




\subsection{Main results}
\label{MR}
A main question in interface models is
whether there exists (maybe under some additional assumptions on the potential $V$ and on the Gibbs measure)
a unique infinite-volume Gibbs measure (or gradient Gibbs measure) describing a
localized interface. 

When there is no disorder, it is known that the Gibbs measure $\nu[\xi=0]$ does not exist in infinite-volume for $d=1,2$,
but the gradient Gibbs measure $\mu[\xi=0]$ does exist in infinite-volume for $d\ge 1$. 
Regarding the uniqueness of gradient Gibbs measures, Funaki and Spohn \cite{FS} showed 
that for uniformly strictly convex potentials $V$ a gradient Gibbs measure $\mu[\xi=0]$ is uniquely determined by
the \textit{tilt} $u\in\R^d$. This result has been extended to a certain class of non-convex potentials by Cotar and Deuschel in \cite{CD}.

For (strongly) non-convex $V$, new phenomena appear: There is
a first-order phase transition from uniqueness to non-uniqueness of the Gibbs measures (at tilt zero), as shown in \cite{BK} and \cite{CD}. More precisely, the model considered in \cite{BK} has potentials of form 
\begin{equation}
\label{Family}
e^{-V_b(\eta_b)}:=p e^{-\k'_b(\eta(b))^2}+(1-p)e^{-{\k''}_b(\eta(b))^2 },~\k'_b,\k''_b>0, p\in [0,1].
\end{equation}
The authors prove in \cite{BK} that there are deterministic choices of $\k'_b, {\k}''_b,p$, independent of the  
bonds $b$, such that there is phase coexistence for the gradient measure 
with tilt $u=0$. 
On the other hand, in \cite{CD} uniqueness is proved for the same potential for different values of $\k', {\k}'',p$ and for $u\in\R^d$. The transition is due to the temperature which changes the structure of the
interface. This
phenomenon is related to the phase transition seen in rotator models with
very nonlinear potentials exhibited in \cite{ESHL1} and \cite{ESHL2}, where the basic mechanism is an
energy-entropy transition. 

How does disorder change these results? In \cite{KO} the authors showed 
that for model A there is no disordered
infinite-volume random Gibbs measure for $d = 1,2$, which is not surprising
since there exists no Gibbs measure without
disorder. Surprising is that, as shown in \cite{EK}, for model A there is also no disordered shift-covariant gradient Gibbs measure when $d=1,2$, and no disordered Gibbs measures for $d=3,4$, as shown in \cite{CK}. For model B, one can reason similarly as for $d=1,2$ in model A (see Theorem 1.1 in \cite{KO}) to show 
that there exists no infinite-volume random Gibbs measure if $d=1,2$. Concerning the question of existence of shift-covariant gradient Gibbs measures, we proved in \cite{CK} that there exists at least one shift-covariant gradient Gibbs measure: for model A when $d\ge 3$ and $\E(\xi(0))=0$, and for model B when $d\ge 1$. 

In this paper, we are interested under what conditions there exists a unique random infinite 
volume gradient Gibbs measure for the two models.

\vspace{1.5mm}
Before we state our main results, we will introduce one more definition.

\begin{defn}
\label{ergodic}
A measure $\P$ is ergodic with respect to translations of $\zd$,  that is $\P\circ (\tau_v)^{-1}=\P$ for all $v\in\zd$ and $\P(A)\in\{0,1\}$ for all $A\in {\cal F}$ such that $\tau_v(A)=A$ for all $v\in\zd$ (for the definition and main theorems of ergodic measures see, for example, Definition 2.3 in \cite{FL} and Chapter 14 in \cite{giorgii}).
\end{defn}

The uniqueness theorem we are about to prove reads as follows. 

\begin{thm} 
\label{uniq}  
Let $u\in\R^d$. 
\begin{enumerate}
\item [(a)] \textbf{(Model A)} Let $d\ge 3$. Assume that $V$ satisfies (\ref{tag22}) and that $(\xi(x))_{x\in\zd}$ have symmetric distributions. For $d=3$ we will also assume that the distribution of $\xi(0)$ satisfies (\ref{gap}). Then there exists a $\P$-almost surely unique shift-covariant gradient Gibbs  measure $\xi\rightarrow \mu^u[\xi]$ defined as in Definition \ref{shiftcov1} with expected tilt $u$, that is with
\begin{equation}
\label{tilt1}
 \E\left(\int\mu^u[\xi](\ormd\eta)\eta(b)\right)=\langle u,y_b-x_b \rangle ~\mbox{for all bonds}~b=(x_b,y_b)\in(\Z^d)^*,
\end{equation}
which satisfies the integrability condition
\begin{equation}
\label{intcond1}
\E\int\mu^u[\xi](\ormd\eta)(\eta(b))^2<\infty~\mbox{for all bonds}~b\in (\zd)^*,
\end{equation}
and such that the annealed measure ${\mu}_{av}^u(\rmd\eta):=\E\int\mu^u[\xi](\ormd\eta)$ is ergodic under the shifts $\{\tau_v\}_{v\in\zd}$. 

\item[(b)] \textbf{(Model B)} Let $d\ge 1$. Assume that for $\P$-almost every $\omega$, $V^\omega_{(x,y)}$ satisfies (\ref{tag23}) uniformly in the bonds $(x,y)$. Then 
there exists a $\P$-almost surely unique shift-covariant gradient Gibbs  measure $\omega\rightarrow\mu^u[\omega]$ defined as in Definition \ref{shiftcov1} with expected tilt $u$, that is with
\begin{equation}
\label{tilt2}
 \E\left(\int\mu^u[\omega](\ormd\eta)\eta(b)\right)=\langle u,y_b-x_b \rangle ~\mbox{for all bonds}~b=(x_b,y_b)\in(\Z^d)^*,
\end{equation}
which satisfies the integrability condition
\begin{equation}
\label{intcond2}
\E\int\mu^u[\omega](\ormd\eta)(\eta(b))^2<\infty~\mbox{for all bonds}~b\in (\zd)^*,
\end{equation}
and such that the annealed measure ${\mu}_{av}^u(\rmd\eta):=\E\int\mu^u[\omega](\ormd\eta)$ is ergodic under the shifts $\{\tau_v\}_{v\in\zd}$.
\end{enumerate}
In words, uniqueness holds for both models in the class of shift-covariant gradient Gibbs measures with ergodic annealed measure and given expected tilt $u$, which class is shown to be non-empty.
\end{thm}
Before we proceed, we note the following
\begin{rem}
\begin{itemize}
\item [(a)] Condition (\ref{tilt1}) (respectively (\ref{tilt2})) is logically stronger than saying that \\
``$\mu[\xi](\nabla_i\phi(x))= \mu'[\xi](\nabla_i\phi(x))$, for all $x\in\zd$, $i\in\{1,2,\ldots, d\}$, 
and  for $\P$-almost every $\xi$, implies that $\mu[\xi]=\mu'[\xi]$". The latter statement would just say 
that the one-dimensional random marginals of the disorder-dependent gradient Gibbs measure $\xi\rightarrow\mu[\xi]$
determine the measure, our theorem 
says that an average tilt determines the measure. 
\item [(b)]  Consider on the other hand a disordered model corresponding to the (very) non-convex potential in (\ref{Family}). Choose $\k'_b$ and/or ${\k}''_b$ random with bounded support, bounded against $0$ from below. 
We may just make one of them random, say $\k'_b$ for instance, or take
$\k'_b=\k' + \o_b$, ${\k}''_b={\k}''+ \o_b$, with $\o_b$ random. 
Then
\begin{equation}
\label{Family1}
e^{-V_b(\eta(b))}:=e^{-\o_b(\eta(b))^2}( 
p e^{-\k'(\eta(b))^2}+(1-p)e^{-{\k}''(\eta(b))^2 }).  
\end{equation}
According to Theorem 3.1 and Remark \ref{rema} c) below, we have existence of a shift-covariant random gradient 
measure with given direction-averaged tilt. Then intuitively one could think that 
an adaptation of the Aizenman-Wehr argument in \cite{AW} (which poses serious
problems in our case because of the unboundedness of the perturbation 
$e^{-\o_b(\eta(b))^2}$) 
should say that when there are two hypothetical gradient 
measures $\mu(\o)$ and $\bar{\mu}(\o)$ with equal expected value 
$\E\mu(\eta(b))=\E\bar\mu(\eta(b))$, the measures are the same in low dimensions, unlike for the equivalent model without disorder,
while one could imagine that in sufficiently high dimensions they are different. 

\end{itemize}
\end{rem}
The deduction of Theorem \ref{uniq} relies partly on a subtle modification of the method of Funaki and Spohn for gradients without disorder from Theorem 2.1 in \cite{FS}, and differs significantly in two main aspects from the proof therein. More precisely, we are able to use neither the shift-invariance and ergodicity of the disordered gradient Gibbs measures nor the extremal/ergodic decomposition of shift-invariant Gibbs measures, which are two main ingredients used in the proof of Theorem 2.1 in \cite{FS}, as in our case the random gradient Gibbs measures are neither ergodic, nor shift-invariant. Furthermore, we are unable to use arguments similar to the ones in \cite{FS} - used there for the case without disorder to construct an ergodic gradient Gibbs measure. It is also worth mentioning here that we cannot assume a priori that there exists a random gradient Gibbs measure - with or without given expected tilt - which is $\P$-a.s. extremal, or which has the property that the corresponding averaged-over-the-disorder measure is ergodic. 
It seems difficult to construct a $\P$-a.s. extremal random gradient Gibbs measure; for example, since the FKG inequality fails in uniformly strictly convex regime for the finite-volume gradient Gibbs measure, we lack monotonicity arguments as used, for example, for the random-field Ising model in Corollary 4.3 from \cite{AW} for such a construction. 
Moreover, the lack of shift-invariance of the disordered gradient Gibbs measure causes serious complications for the arguments necessary to prove Theorem \ref{uniq}. 

One of the main ingredients in our proof is Theorem \ref{directilt}, a far from trivial result of a.s. existence of a shift-covariant gradient Gibbs measure with given direction-averaged tilt, proved by means of the Brascamp-Lieb inequality and (for model A) also of a Poincar\'e-type inequality. We will then exploit in Lemma \ref{annealeduniq} the rapid decay of the norm $\|\eta\|_r, r>0$, and use Theorem \ref{directilt}, to obtain uniqueness of the averaged-over-the-disorder gradient Gibbs measure (the \textit{annealed} measure) with given \textit{direction-averaged tilt}. Together with Proposition \ref{kom} - which is the key to allowing us to pass from uniqueness of the annealed measure to almost sure uniqueness of the corresponding disorder-dependent, gradient Gibbs measure (the \textit{quenched} measure) - Lemma \ref{annealeduniq} will provide us with the statement from Theorem \ref{extremeuniq}, of uniqueness of the \textit{quenched} gradient Gibbs measure with given \textit{direction-averaged expected tilt}. From this last theorem we will also derive the ergodicity of the annealed gradient Gibbs measure with given direction-averaged tilt. We will then upgrade the result in Theorem \ref{extremeuniq} to the statement from Theorem \ref{uniq} of uniqueness with given \textit{expected tilt} and corresponding \textit{ergodic} annealed measure. 
\vspace{2mm}

Let $C^1_b(\chi_r)$ denote the set of differentiable functions depending on finitely many 
coordinates with bounded derivatives, where $\chi_r$ was defined in Subsection 1.2.2. Let $F\in C^1_b(\chi_r)$.
We denote by
\begin{equation}
\label{partialbond}
\partial_bF(\eta):=\frac{\partial F(\eta)}{\partial \eta(b)}~\mbox{and}~||\partial_bF||_{\infty}:=\sup_{\eta\in\chi}|\partial_bF(\eta)|.
\end{equation}
Let $b=(x_b,y_b)\in (\zd)^*$. In the formulas below, and to avoid exceptional cases when $b=0$, we denote by $]|b|[=\max\{|x_b|,1\}$, where $|x_b|$ is the Euclidian norm. We prove next the decay of covariance \textit{with respect to the averaged-over-the-disorder} random gradient Gibbs measure from Theorem \ref{uniq}.
\begin{thm}
\label{decay}
Let $u\in\R^d$. 
\begin{enumerate}
\item [(a)] \textbf{(Model A)} Let $d\ge 3$. Assume that $V$ satisfies (\ref{tag22}) and that $(\xi(x))_{x\in\zd}$ are i.i.d with mean $0$ and the distribution of $\xi(0)$ satisfies (\ref{gap}). Then if $\xi\rightarrow\mu^u[\xi]$ is any shift-covariant gradient Gibbs measure constructed as in \cite{CK}, $\xi\rightarrow\mu^u[\xi]$ satisfies the following decay of covariances for all $F,G\in C_b^1(\chi_r)$
$$|\Cov\left(\mu^u[\xi](F(\eta)),\mu^u[\xi](G(\eta))\right)|\le c\sum_{b,b'\in (\zd)^*}\frac{||\partial_b F||_\infty||\partial_{b'} G||_\infty}{]|b-b'|[^{d-2}},$$
for some $c>0$ which depends only on $d, C_1$, $C_2$ and on the number of terms $b,b'$ in $F$ and $G$.
\item[(b)] \textbf{(Model B)} Let $d\ge 1$. Even though we can consider more general disorder structures, we assume for simplicity that $V^\omega_{(x,y)}(\phi(x)-\phi(y))=V_{(x,y)}(\omega(x,y),\phi(x)-\phi(y))$ and that for all $b=(x,y)\in (\zd)^*$ there exists 
$\frac{\partial^2V^\omega_{(x,y)}}{\partial\omega(b)\eta(b)}$ with $\left|\frac{\partial^2V^\omega_{(x,y)}}{\partial\omega(b)\eta(b)}\right|\le f_{1,b}(\omega)\left|\eta(b)\right|+f_{2,b}(\omega)$ for some measurable $f_{i,b}:\Omega\rightarrow\R_{+}$ with $\sup_b\E(f_{i,b}^p)<\infty, 2<p<\infty, i=1,2$. 
Assume also that $\omega(x,y)$ are i.i.d. for all $(x,y)$, that the distribution of $\omega(x,y)$ satisfies (\ref{gap}) and that $V^\omega_{(x,y)}$ satisfies (\ref{tag23}) for $\P$-almost every $\omega$ and uniformly in the bonds $(x,y)$. Then if $\omega \rightarrow\mu^u[\omega]$ is any shift-covariant gradient Gibbs measure constructed as in \cite{CK} ($\P$-almost surely unique by Theorem \ref{uniq}), $\omega\rightarrow\mu^u[\omega]$ satisfies the following decay of covariances for all $F,G\in C_b^1(\chi_r)$ 
$$| \Cov\left(\mu^u[\omega](F(\eta)),\mu^u[\omega](G(\eta))\right)|\le c\sum_{b,b'\in (\zd)^*}\frac{||\partial_b F||_\infty||\partial_{b'} G||_\infty}{]|b-b'|[^{d}},$$
for some $c>0$ which depends only on $d, C_1$, $C_2$  and on the number of terms $b,b'$ in $F$ and $G$.
\end{enumerate}
\end{thm}

\begin{rem}
We note here that one can easily verify in the case with quadratic potentials that the above bounds are optimal by simple Gaussian computations. Moreover, for model A
one can prove the following for $F=G=V'$ and for large enough $|b-b'|$, by generalizing the proof of Theorem 1.2 in \cite{EK} from $d=3$ 
to any dimension $d\geq 3$: An upper bound of form
\begin{equation}
\label{contrabound}\begin{split}
&|\Cov\left(\mu^u[\xi](V'(\eta(b))),\mu^u[\xi](V'(\eta(b')))\right)|\ \leq \text{Const} \,\, ]|b-b' |[^{-q},~~q>0,\cr
\end{split}
\end{equation}
cannot be  
true for q  $\geq d-2$.
In words, there cannot be a uniform upper bound with a better exponent. However, this does not exclude that some of the covariances for specifically  
chosen bonds $b, b'$ might even be zero. The statement holds even for \textbf{highly non-convex} potentials like the one in \cite{BK}.

To prove this, we assume an upper bound $q$ and we will show that it cannot be greater than $q=d-2$. The proof follows from the identity (18)  in \cite{EK}. 
This identity is obtained from a spatial sum of the divergence equation (15), it holds for arbitrary volumes, and is independent of the spatial dimension. 
Considering balls of radius $L$ one derives that, for $L$ large enough, the assumed decay would imply $L^d\leq \bar{c} L^{2(d-1)-q}$, for some $\bar{c}>0$ depending on $d$, 
which proves the desired bound on $q$.
\end{rem}
\begin{rem}
In view of \cite{mar} and of \cite{men}, it would be possible to weaken the i.i.d. assumption on the disorder from Theorem \ref{decay} to certain weak dependence and stationarity assumptions. However, for simplicity of calculations purposes, we will restrict ourselves to the i.i.d. case.
\end{rem}


The methods we employ for our main theorems can be used to tackle similar questions for other gradient models with disorder such as, for example, the gradient model on the supercritical percolation cluster from \cite{CY} or the gradient model with disordered pinning from \cite{CM}.

\vspace{2mm}

The rest of the paper is organized as follows: In Section 2 we recall a number of basic definitions and main properties used in the proof of our main results. In Section 3, we show in Theorem \ref{directilt}  one of the main ingredients necessary for the proof of Theorem \ref{uniq}, the existence of a shift-covariant gradient Gibbs measure with given direction-averaged tilt. In Section 4, we upgrade in Theorem \ref{extremeuniq} this statement of existence to one of uniqueness of measures with given direction-averaged tilt, which implies also the ergodicity of the corresponding annealed measure in Theorem \ref{ergodic}. In Section 5, we prove the decay of covariances result from Theorem \ref{decay}.

\section{Preliminary notions}
For the reader's convenience, we will introduce in this section a number of notions and results used in the proofs of our main statements, Theorems \ref{uniq} and \ref{decay}.

\subsection{Estimates for the discrete Green's functions on $\Z^d$}

\label{greenfunction}
We will state first a probabilistic interpretation of the discrete Green's function. Let $A$ be an arbitrary subset in $\Z^d$ and let $x\in A$ be fixed. Let $\P_x$ and $\E_x$ be the probability law and expectation, respectively, of a simple random walk $X:=(X_k)_{k\geq 0}$ starting from $x\in\Z^d$; the \textit{discrete Green's function $G_A(x,y)$} is the expected number of visits to $y\in A$ of the walk $X$ killed as it exits $A$, i.e.
$$G_A(x,y)=\E_x\left[\sum_{k=0}^{\tau_A-1} 1_{(X_k=y)}\right]=\sum_{k=0}^\infty\P_x(X_k=y,k<\tau_A),~~y\in\Z^d,$$ 
where $\tau_A=\inf\{k\ge 0:X_k\in A^c\}$ is the first exit time of $X_k$ from $A$. 

We will next give  some well-known properties of the Green's functions. To avoid exceptional cases when $x=0$, let us denote by $]|x|[=\max\{|x|,1\}$, where $|x|$ is the Euclidian norm. Let $\Lambda_N=[-N,N]^d$.
\begin{prop}
\label{properties}
\begin{enumerate}
\item [(i)] If $d\ge 3$, then $\lim_{N\rightarrow\infty}G_{\Lambda_N}(x,y):=G(x,y)$ exists for all $x,y\in\Z^d$ and as $|x-y|\rightarrow\infty$,
$$G(x,y)=\frac{a_d}{|x-y|^{d-2}}+O(|x-y|^{1-d}),$$
with $a_d=\frac{2}{(d-2)w_d}$, where $w_d$ is the volume of the unit ball in $\R^d$.
\item [(ii)] Let $B_r=\{x\in\Z^d:|x|<r\}$; then for $x\in B_N$
$$G_{B_N}(0,x)=
\left\{
\begin{array}{lcc}
\frac{2}{\pi}\log\frac{N}{]|x|[}+o\left(\frac{1}{]|x|[}\right)+O\left(\frac{1}{N}\right)& \mbox{if} & d=2\\
\frac{2}{(d-2)w_d}\left[~]|x|[^{2-d}-N^{2-d}+O\left(~]|x|[^{1-d}\right)\right] & \mbox{if} & d\ge 3.
\end{array}
\right.$$
Let $\epsilon>0$. If $x\in B_{(1-\epsilon)N}$ the following inequalities hold:
$$G_{B_{{\epsilon N}}}(0,0)\le G_{B_N}(x,x)\le G_{B_{2N}}(0,0).$$
\item [(iii)] $G_A(x,y)=G_A(y,x)$.
\item [(iv)]  $G_A(x,y)\le G_B(x,y)$, if $A\subset B$.
\end{enumerate}
\end{prop}
\noindent For proofs of (i), (iii) and (iv) from Proposition \ref{properties} above we refer to Chapter 1 from \cite{Lawl1} and for proof of (ii) we refer to Lemma 1 from \cite{L}.

\subsection{Covariance inequalities}

We will state next some variance and covariance inequalities  for finite-volume Gibbs measures, needed for the proof of our main results Theorem \ref{uniq} and Theorem \ref{decay}. Following \cite{DGI}, we will state these inequalities for the Hamiltonian
\begin{equation}
\label{genhamil}
H_{\Lambda}^\psi(\phi)[\vartheta]:=\frac{1}{2}\sum_{x,y\in\Lambda,|x-y|=1}V_{(x,y)}(\phi(x)-\phi(y))+\sum_{x\in\Lambda,y\in\partial\Lambda,|x-y|=1}V_{(x,y)}(\phi(x)-\psi(y))+\sum_{x\in\Lambda}\vartheta(x)\phi(x),
\end{equation}
which, for \textit{fixed} disorder, covers \textit{both} the cases of our models (A) and (B). We assume that the external field $(\vartheta(x))_{x\in\zd}\in\R^{\zd}$. We have the usual conditions on $V_{(x,y)}$: for some given $0<C_1<C_2$, $V_{(x,y)}$ obey the following bounds, uniformly in the bonds $(x,y)$
\begin{equation}
\label{ellipt}
C_1\le (V_{(x,y)})''(s)\leq C_2 ~\mbox{for all}~s\in\RR.
\end{equation}
We assume also that for each bond $(x,y)$,  $V_{(x,y)}\in C^2(\Bbb R)$ is an even function. We define $\nu^\psi_\Lambda[\vartheta]$ and $\mu^\rho_\Lambda[\vartheta]$ corresponding to $H_{\Lambda}^\psi(\phi)[\vartheta]$ as in Subsection \ref{AG}.

\subsubsection{Helffer-Sj\"ostrand (random walk) representation}
\label{rwrep0}

The idea, due to Helffer-Sj\"ostrand, originally developed in [15] and reworked probabilistically
in \cite{DGI}, \cite{gos}, is to describe the correlation functions under the
Gibbs measures in terms of the first exit distribution and occupation time of a certain random walk in random environments. More precisely, given the \textit{time-independent} environment $\{\nabla\phi\}$, we will denote by $\{X_t,t\ge 0\}$ the random walk on $\zd$ with time-dependent jump rates along the bond $b = (x_b, y_b)\in (\zd)^*$ given by 
$$a^{\nabla\phi}(t,x_b,y_b)=V''_b(\phi_t(x_b)-\phi_t(y_b)).$$

Since the function V is even, we have symmetric jump rates: $a^{\nabla\phi}(t,x_b,y_b)=a^{\nabla\phi}(t,y_b,x_b)$. Moreover the condition (\ref{ellipt}) guarantees ellipticity, so our random walk exists. We write next the transition probability of the random walk killed at the time when it goes outside of $\Lambda$
\begin{equation}
\label{dgirep1a}
p_\Lambda^{\nabla\phi}(s,x,t,y):=\mathbb{P}^{\nabla\phi}(X_t=y,~t<\tau_\Lambda|X_s=x)~~\mbox{and}~~g^{\nabla\phi}_\Lambda(x,y)=\int_0^\infty p_\Lambda^{\nabla\phi}(0,x,t,y)\rmd t,
\end{equation}
where, as before, $\tau_\Lambda:=\inf\{i>0,X_i\in\Lambda^c\}$ and $t\ge s\ge 0$. We note here that $p_\Lambda^{\nabla\phi}(s,x,t,y)$ depends on $\nabla\phi$ \textit{only} through $a^{\nabla\phi}$. We now have from Proposition 2.2 in \cite{DGI} (see also Theorem 4.2 in \cite{FL})
\begin{prop}\textbf{(Random walk representation)}
\label{dgirep}
Fix $\Lambda\subset\zd$ finite and $\psi\in \R^{\zd}$. Let $F,G$ be the set of differentiable functions on $\R^{\Lambda}$ with bounded derivatives. Then
\begin{equation}
\label{reprep}
\cov_{\nu^\psi_\Lambda[\vartheta]}(F(\phi),G(\phi))=\int_{0}^\infty\sum_{x,y\in\Lambda}\E_{\nu^\psi_\Lambda[\vartheta]}\left(\partial_x F(\phi)\partial_y G(\phi)p_\Lambda^{\nabla\phi}(0,x,s,y)\right)\rmd s,
\end{equation}
where we denoted by $\partial_xF(\phi):=\frac{\partial F(\phi)}{\partial\phi(x)}$, and by $\E_{\nu^\psi_\Lambda[\vartheta]}$ and $\cov_{\nu^\psi_\Lambda[\vartheta]}$ the expectation, respectively covariance, with respect to $\nu^\psi_\Lambda[\vartheta]$. 
In the special case that $F(\phi)=\phi(a)$ and $G(\phi)=\phi(b)$ for some $a,b\in\Lambda$, we simply have
\begin{equation}
\label{covrwrep}
\cov_{\nu^\psi_\Lambda[\vartheta]}(\phi(a), \phi(b))=\int_{0}^\infty\E_{\nu^\psi_\Lambda[\vartheta]}\left(p_\Lambda^{\nabla\phi}(0,a,s,b)\right)\rmd s\le \int_{0}^\infty\E_{\nu^\psi_\Lambda[\vartheta]}\left(p^{\nabla\phi}(0,a,s,b)\right)\rmd s.
\end{equation}
\end{prop}
Let us now define
\begin{equation}
\label{dgirep1}
p^{\nabla\phi}(s,x,t,y):=\lim_{|\Lambda|\rightarrow\infty}p_\Lambda^{\nabla\phi}(s,x,t,y)=\mathbb{P}^{\nabla\phi}(X_t=y|X_s=x)~~\mbox{and}~~g^{\nabla\phi}(x,y)=\int_0^\infty p^{\nabla\phi}(0,x,t,y)\rmd t.
\end{equation}
We note here that in the case with $\vartheta=0$, there exists for all $u\in\R^d$ a \textit{unique shift-invariant extremal} infinite-volume gradient Gibbs measure $\mu^u[\vartheta=0]$ with tilt $u$ (as proved in \cite{FS}), which satisfies a random walk representation as in Proposition \ref{dgirep} above, with $p^{\nabla\phi}$ replacing $p^{\nabla\phi}_\Lambda$ in (\ref{reprep}) (for a statement see, for example, Proposition 3.1 in \cite{gos} or (6.7) in \cite{DD}). However, the extension to infinite-volume is non-trivial and, unlike the corresponding finite-volume representation, the proofs rely on the extremality of $\mu^u[\vartheta=0]$. 

\vspace{1.5mm}

We will use in our proof of Theorem \ref{directilt} (a) and Theorem \ref{decay} the following properties of  $g_\Lambda^{\nabla\phi}(x,z)$ and $g^{\nabla\phi}(x,z)$, well-known in the gradient literature and stated here for the reader's convenience.
\begin{prop}
\label{propg}
 Let $d\ge 3$. 
\begin{enumerate}
\item [(i)] There exist $c_{-},c_{+}>0$, which depend only on $d, C_1$ and $C_2$, such that for all $x,z\in\zd$, $\nabla\phi\in (\zd)^*$ and $\Lambda\subset\zd$ finite, we have
$$0\le g_\Lambda^{\nabla\phi}(x,z)\le \frac{c_{+}}{]|x-z|[^{d-2}}~~~~~~\mbox{and}~~~~~~\frac{c_{-}}{]|x-z|[^{d-2}}\le g^{\nabla\phi}(x,z)\le \frac{c_{+}}{]|x-z|[^{d-2}}.$$
\item [(ii)] There exists  $c_{+}>0$, which depends only on $d, C_1$ and $C_2$, such that for all $x,z\in\zd$, $\rho\in (\zd)^*$ and $\Lambda\subset\zd$ finite, we have
$$0\le \cov_{\nu^\psi_\Lambda[\vartheta]}(\phi(x),\phi(z))\le \frac{c_{+}}{]|x-z|[^{d-2}}.$$
\item [(iii)]  There exist $\tilde{C}(d),\rho>0$, which depends only on $d, C_1$ and $C_2$, such that for all $R>0,\Lambda\subset\zd$ finite, $\nabla\phi\in (\zd)^*$, $z\in\zd$ and all $\alpha,\beta\in \{1,2,\ldots, d\}$, we have
\begin{equation}
\label{gogo30a}
\sum_{x:R\le |x-z|\le 2R}\left(g_\Lambda^{\nabla\phi}(x, z)-g_\Lambda^{\nabla\phi}(x+e_\alpha, z)\right)^2\le \tilde{C}(d)R^{2-d},
\end{equation}
and (for $d\ge 1$)
\begin{eqnarray}
\label{gogo30}
\sum_{x:R\le |x-z|\le 2R}\left(g^{\nabla\phi}_\L(x, z)-g^{\nabla\phi}_\L(x+e_\alpha, z)-g^{\nabla\phi}_\L(x, z+e_\beta)+g^{\nabla\phi}_\L(x+e_\alpha, z+e_\beta)\right)^2
\le \tilde{C}(d)R^{-\rho},\nonumber\\
\end{eqnarray}
where $e_\alpha$ and $e_\beta$ are the unit vectors in direction $\alpha$, respectively $\beta$. Note that (\ref{gogo30}) can be proved in a stronger form for $d\ge 2$ (i.e., with the suboptimal bound $R^{2-d-\rho}$).  
\item [(iv)] There exist $\delta,{C_+}>0$, which depend only on $d, C_1$ and $C_2$, such that for all $\Lambda\subset\zd$ finite, $\nabla\phi\in (\zd)^*$, $z\in\zd$ and all $\alpha\in \{1,2,\ldots, d\}$, we have
\begin{equation}
\label{nash}
\left|g_\Lambda^{\nabla\phi}(x, z)-g_\Lambda^{\nabla\phi}(x+e_\alpha, z)\right|\le \frac{C_{+}}{]|x-z|[^{d-2+\delta}}.
\end{equation}
\item [(v)] Let $\gamma$ be a shift-invariant measure on $\chi$, let $d\ge 1$ and let $1\le p<\infty$. There exists $\bar{C}>0$, which depends only on $d,p, C_1$ and $C_2$, such that for all $\Lambda\subset\zd$ finite, $\nabla\phi\in (\zd)^*$, $z\in\zd$ and for all $\alpha,\beta\in \{1,2,\ldots, d\}$, we have 
\begin{equation}
\label{ahr}
\gamma\left(\left(g^{\nabla\phi}(x, z)-g^{\nabla\phi}(x+e_\alpha, z)\right)^{2p}\right)\le \frac{\bar{C}}{]|x-z|[^{2pd-2p}}.
\end{equation}
and
\begin{eqnarray}
\label{ahr1}
\gamma\left(\left(g^{\nabla\phi}(x, z)-g^{\nabla\phi}(x+e_\alpha, z)-g^{\nabla\phi}(x, z+e_\beta)+g^{\nabla\phi}(x+e_\alpha, z+e_\beta)\right)^{2p}\right)\le \frac{\bar{C}(d)}{]|x-z|[^{2pd}}.\nonumber\\
\end{eqnarray}
\end{enumerate}
\end{prop}

{\bf Proof.}

For a proof of (i), (and in view of the classical De Giorgi-Nash-Moser theory), see for example Propositions B.3 and B.4 in \cite{gos}. To prove (ii), we combine (\ref{covrwrep}) from Proposition \ref{dgirep} with Proposition \ref{propg} (i) (see Theorem 4.13 in \cite{FL} for an extended proof of (ii)). The proof of (\ref{gogo30a}) in (iii) relies on a standard Caccioppoli argument with respect to $x$, and is based on the decay of $g_\Lambda^{\nabla\phi}(x+e_\alpha, z)$ given in (i) (for a similar proof and discussion, see for example Lemma 2.9 in \cite{GO1}; for a statement of Caccioppoli's inequality, see for example Propositions 2.1 and 4.1 in \cite{DD}). For a proof of (\ref{gogo30}), see (30) in Lemma 6 from \cite{MAO}. The stronger form of (\ref{gogo30}) for $d\ge 2$ (i.e., with the suboptimal bound $R^{2-d-\rho}$) can be proved by means of (\ref{gogo30}) and of Caccioppoli's inequality (see the explanation in Section 7.2 from \cite{MAO}). The proof of (iv) follows from the famous Nash continuity estimate, as stated for example in Proposition B.6 from \cite{gos}. For a proof of (v), see Theorem 1 from \cite{MAO}.

See also \cite{GO1} and \cite{MAO} for more estimates and extended explanations on $p^{\nabla\phi}(0,x,z)$ and $g^{\nabla\phi}(x,z)$.

$\Cox$.


\subsubsection{The Brascamp-Lieb inequality}

The Brascamp-Lieb inequality states that for $\g$ a centered Gaussian distribution on $\R^N, N\ge 1$,
and $\mu$ a distribution on $\R^N$ such that there exists $d\mu/d\gamma = e^{-f}$ for a convex function $f$, one has for all $v\in\RR^N$ and for all convex real functions 
$L$, bounded below, that 
\begin{equation}
\label{111}
\mu \Big(L\bigl(v \cdot (X -\mu(X))\bigr)\Big)\leq \g\Big(L(v \cdot X)\Big).
\end{equation}
The above is the formulation by Funaki in \cite{FL}. An application of (\ref{111}) to our $\mu^\rho_\Lambda[\vartheta]$ case with $L(s)=s^2$ (see also Lemma 2.8 in \cite{DGI} for the proof in the case with $f$ equal to $H_\Lambda^\psi[\vartheta]$ as in (\ref{genhamil})), would give for example that 
\begin{equation}
\label{bl11}
\mu^{\rho}_{\L}[\vartheta] \Bigl( \Bigl[
\phi(x_0)-\phi(y_0)- \mu^{\rho_u}_{\L}[\x] \bigl( \phi(x_0)-\phi(y_0)\bigr)\Bigr]^2\Bigr)
\le \frac{1}{C_1} \mu_{G,{\L}}^{\rho}[\vartheta=0]\Bigl(  \Bigl[\phi(x_0)-\phi(y_0)\Bigr]^2\Bigr),
\end{equation}
where $\mu^{\rho}_{G, {\L}}[\vartheta=0]$ is the corresponding Gaussian gradient Gibbs measure with potential $V_0(s)=\frac{s^2}{2}$ and external field $\vartheta=0$.

\subsubsection{Localization of the variance under pinning}

A crucial property of low-dimensional ($d = 1, 2$)
continuous interfaces \textit{without disorder} is that the local variance of the field has a slow growth. However, it turns out that pinning a single point is sufficient to \textit{localize} the
field, in the sense that an infinite-volume Gibbs measure exists. More precisely,
let us consider the Gaussian measure $\nu^{0}_{G, {\L_N}\setminus \{0\}}[\vartheta=0]$, i.e. the Gaussian Gibbs measure with $0$ boundary conditions outside $\Lambda_N:=[-N,N]^d$ and at the origin. Then one can show that for any $a\in\zd$, we have
$$\lim_{N\rightarrow\infty}\var_{\nu^{0}_{G, {\L_N}\setminus \{0\}}[\vartheta=0]}(\phi_a)\simeq |a|~~\mbox{if}~~d=1~~~\mbox{and}~~~\lim_{N\rightarrow\infty}\var_{\nu^{0}_{G, {\L_N}\setminus \{0\}}[\vartheta=0]}(\phi_a)\simeq \log|a|~~\mbox{if}~~d=2.$$
Actually, one even has that
\begin{equation}
\label{varloc}
\sup_{a\neq 0}\frac{\lim_{N\rightarrow\infty}\var_{\nu^{0}_{G, {\L_N}\setminus \{0\}}[\vartheta=0]}(\phi_a)}{\var_{\nu^{0}_{G, {\L(a)}}[\vartheta=0]}(\phi_a)}\simeq 1,
\end{equation}
where $\Lambda(a)=\{b\in\zd:|a-b|_\infty|\le |b|_\infty\}$. In the above, $\simeq$ stands for a multiplicative constant which only
depends on the dimension $d$.

In the above, we have taken 0 boundary conditions outside $\Lambda_N$, but any boundary
conditions not growing too fast with $N$ would have given the same result. For more on the above estimates and localization of the variance under pinning in general, see for example \cite{velenik}.

\subsection{Covariance inequalities under the disorder}

Similarly to the proof of Lemma 3 from \cite{GO1} we have the following covariance inequality,  which in the particular case of the variance is a weakened version of a second order Poincar\'e inequality.
\begin{prop}
\label{golemma}
Fix $n\in\N$ and let $a = (a_i)_{i=1}^n$ be a sequence of independent random variables with uniformly-bounded finite second moments on a probability space $(\Omega, {\cal F},\P)$. Let $X,Y$ be Borel measurable functions of $a\in\R^n$ (i.e. measurable w.r.t.
the smallest $\sigma$-algebra on $\R^N$ for which all coordinate functions $\R^n\ni a\rightarrow a_i\in\R$ are Borel
measurable). Then we have
\begin{equation}
\label{poincarre}
\left|\cov(X,Y)\right|\le \max_{1\le i\le n} \var(a_i)\sum_{i=1}^n\left(\int \sup_{a_i}\left|\frac{\partial X}{\partial a_i}\right|^2 \rmd\P\right)^{1/2} \left(\int \sup_{a_i}\left|\frac{\partial Y}{\partial a_i}\right|^2 \rmd\P\right)^{1/2},
\end{equation}
where $\sup_{a_i} \left|\frac{\partial Z}{\partial a_i}\right|$ denotes the supremum of the modulus of the $i$-th partial derivative
$$\frac{\partial Z}{\partial a_i}(a_1,\ldots,a_{i-1},a_i,a_{i+1},\ldots, a_n)$$
of $Z$ with respect to the variable $a_i$, for $Z=X,Y$.
\end{prop}
For i.i.d  random variables, one can obtain under the mild assumption (\ref{gap}) on the distribution $\gamma$ of $a_i$ the following stronger variance estimate  
\begin{equation}
\label{ledoux}
\var(X)\le C(d)\sum_{i=1}^n\int \left|\frac{\partial X}{\partial a_i}\right|^2 d\P,
\end{equation}
where $C(d)>0$ depends only on $d$ and on the distribution of $a_i$. 
For the proof of (\ref{ledoux}), see for instance Lemma 1.1 from \cite{Led}; for a related weak dependence statement for absolutely continuous measures, see Theorem 1 from \cite{mar}, for the statement for discrete measures, Theorem 2.1 from \cite{men}.

\subsection{Construction of a shift-covariant random gradient Gibbs measure}

We recall in this subsection the construction of an infinite-volume shift-covariant gradient Gibbs measure,  as given in Theorem 1.7 and in Proposition 3.8 from \cite{CK}. 

Let $u\in\RR^d$ and let the boundary condition $\psi_u(x):=u\cdot x,x\in\zd$. Take $\rho_u(b):=\nabla\psi_u(b)$ for all $b\in (\zd)^*$ and consider the corresponding gradient Gibbs measure $\mu_\L^{\rho_u}[\xi]$ as given by (\ref{fingradgibbs}). Let us now define the \textit{spatially-averaged} measure $\bar\mu^{u}_{\L}[\xi]$ on gradient configurations given by 
\begin{equation}
\label{heldfixed}
\bar\mu^{u}_{\L}[\xi]:=\frac{1}{|\L|}\sum_{x\in \L}\mu^{\rho_u}_{\L+x}[\xi],
\end{equation}
where we defined $\L+x:=\{z+x:z\in\L\}$. This is an extension to our disorder-dependent case 
of the construction of Gibbs measures with symmetries given in \cite{giorgii}, in formula (5.20) from Chapter $5.2$; the construction in \cite{giorgii} was used there to obtain shift-invariant Gibbs measures. We note that in (\ref{heldfixed}), the random field variables $\x$ are held fixed while the volumes $\L+x$ 
are shifted around. From Theorem 1.7 and Proposition 3.8 in \cite{CK}  we have
\begin{prop}\textbf{(Existence of shift-covariant random gradient Gibbs measures)}
\label{existgibbs}
\begin{enumerate}
\item [(a)] \textbf{(Model A)} Let $d\ge 3$ and $\E(\xi(0))=0$. Assume that $V$ satisfies (\ref{tag22}). Then there exists a deterministic subsequence $(m_i)_{i\in\N}$ such that for $\P$-almost every $\xi$
\begin{equation}
\label{hat1}
\hat\mu^{u}_{k}[\xi]:=\frac{1}{k}\sum_{i=1}^k {\bar\mu}^{u}_{\L_{m_{i}}}[\xi]
\end{equation}
converges as $k\to \infty$ weakly to $\mu^u[\xi]$, which is a shift-covariant random gradient Gibbs measure defined as in Definition \ref{shiftcov1}. Moreover, $\mu^u[\xi]$ satisfies the integrability condition
\begin{equation}
\label{intcond}
\E\int\mu^u[\xi](\ormd\eta)(\eta(b))^2<\infty~\mbox{for all bonds}~b\in (\zd)^*.
\end{equation}
\item[(b)] \textbf{(Model B)} Let $d\ge 1$. Assume that for $\P$-almost every $\omega$, $V^\omega_{(x,y)}$ satisfies (\ref{tag23}), uniformly in the bonds. Then there exists a deterministic subsequence $(m_i)_{i\in\N}$ such that for $\P$-almost every $\omega$
\begin{equation}
\label{hat2}
\hat\mu^{u}_{k}[\omega]:=\frac{1}{k}\sum_{i=1}^k {\bar\mu}^{u}_{\L_{m_{i}}}[\omega]
\end{equation}
converges as $k\to \infty$ weakly to $\mu^u[\omega]$, which is a shift-covariant random gradient Gibbs measure defined as in Definition \ref{shiftcov1}
Moreover, $\mu^u[\omega]$ satisfies the integrability condition
\begin{equation}
\E\int\mu^u[\omega](\ormd\eta)(\eta(b))^2<\infty~\mbox{for all bonds}~b\in (\zd)^*.
\end{equation}
\end{enumerate} 
\end{prop}
\begin{rem}
\begin{itemize}
\item [(a)] The above theorem was proved in \cite{CK} without the assumption of strict convexity of the potentials in models (A) and (B).  Note that even though the proofs in \cite{CK} were done under the assumption of i.i.d disorder for both models, only stationarity of the disorder was used in the proofs for model B. Note also that we can also construct the gradient Gibbs measures above  through the use of periodic boundary conditions, which automatically ensures shift-covariance of the quenched measure. 
\item [(b)] Our measures (\ref{hat1}), respectively (\ref{hat2}), are obtained via a construction which resembles the  
construction of the barycenter
of an empirical metastate in the sense of Newman and Stein (see, for example, \cite{SN} for more on this).
The modification we adopted - for the purpose of constructing a shift-covariant random infinite-volume gradient Gibbs measure, as defined in Definition \ref{shiftcov1} - lies in the  
fact that our finite-volume  measures (\ref{heldfixed})
have already undergone a spatial averaging themselves before
they  are summed along the volume sequence indexed by $k$.
\end{itemize}
\end{rem}

\section{Existence of shift-covariant random gradient Gibbs measure with given direction-averaged tilt}

We will prove in this section one of the main ingredients necessary for the proof of our main result in Theorem \ref{uniq}. We will use in our proof the construction of the infinite-volume shift-covariant gradient Gibbs measure from \cite{CK}. 

Fix $u\in\R^d$. We will show that for $\P$-almost every $\xi$ (respectively $\omega$), the following is true: there exists a shift-covariant random gradient Gibbs measure $\mu^u [\xi]$ (respectively $\mu^u[\omega]$), with respect to which the gradient 
averages in any fixed direction $\alpha\in\{1,2,\ldots,d\}$ over the tilt $u$ converge to zero stochastically as $\Lambda\uparrow \Z^d$.  This would exclude that this random gradient Gibbs measure is a linear combination 
between random Gibbs measures which are supported on sets of interfaces 
with two or more different expected tilts. 
More precisely, we will prove
\begin{thm}
\label{directilt}
Fix $u\in\R^d$. 	Let for all $\alpha\in\{1,2,\ldots,d\}$ 
$$E_\alpha:=\{ \eta~|~\lim_{|\Lambda|\rightarrow\infty}\frac{1}{|\Lambda|}\sum_{x\in\Lambda}\eta(b_{x,\alpha})=u_\alpha\},$$
along the sequence of volumes with $b_{x,\alpha}:=(x+ e_\alpha,x)\in (\zd)^*$. 
\begin{enumerate}
\item [(a)] \textbf{(Model A)} Let $d\ge 3$. Assume that $V$ satisfies (\ref{tag22}) and that $(\xi(x))_{x\in\zd}$ have symmetric distribution. For $d=3$ we will also assume that the distribution of $\xi(0)$ satisfies (\ref{gap}). Then there exists a shift-covariant random gradient Gibbs measure defined as in Definition \ref{shiftcov1} which satisfies for $\P$-almost every $\xi$
\begin{equation}
\label{firstlimit}
\mu^u[\xi](E_\alpha)=1,~\alpha\in\{1,2,\ldots,d\}.
\end{equation} 
Moreover, $\mu^u[\xi]$ satisfies the integrability condition
\begin{equation}
\label{intconda}
\E\int\mu^u[\xi](\ormd\eta)(\eta(b))^2<\infty~\mbox{for all bonds}~b\in (\zd)^*.
\end{equation}
\item[(b)] \textbf{(Model B)} Let $d\ge 1$. Assume that for $\P$-almost every $\omega$, $V^\omega_{(x,y)}$ satisfies (\ref{tag23}). Then there exists a shift-covariant random gradient Gibbs measure defined as in Definition \ref{shiftcov1} which satisfies for $\P$-almost every $\omega$
\begin{equation}
\label{firstlimit11}
\mu^u[\omega](E_\alpha)=1,~\alpha\in\{1,2,\ldots,d\}.
\end{equation} 
Moreover, $\mu^u[\omega]$ satisfies the integrability condition
\begin{equation}
\label{intcond2a}
\E\int\mu^u[\omega](\ormd\eta)(\eta(b))^2<\infty~\mbox{for all bonds}~b\in (\zd)^*.
\end{equation}
\end{enumerate}
\end{thm}

{\bf Proof.}

 For both models, we will treat separately in the proof the critical dimensions ($d=3,4$ for model A and $d=1,2$ for model B) where a more delicate analysis is required, and the remaining dimensions. The key idea  to show (\ref{firstlimit}), respectively (\ref{firstlimit11}), is to bound the main quantity to be estimated by a sum of two variances. The first variance can be bound by means of the Brascamp-Lieb inequality and (for $d=1,2$ in model B) also by the variance estimates from (\ref{varloc}). The second variance can be bound for model A by means of Proposition \ref{golemma}; for model B, it will be equal to zero by arguments involving the symmetry of the potentials $V_{(x,y)}$. To further estimate the second variance for model A, we will use the finite-volume random walk representation from Proposition \ref{dgirep}, the bounds from Proposition \ref{propg} (ii), and (for $d=3,4$) also the bounds from Proposition \ref{propg} (iii) and (iv). 
 
By our construction, the tilt $\mu^u[\xi](\ormd\eta)(\eta(b))$ is random for model A, whereas for model B the tilt $\mu^u[\omega](\ormd\eta)(\eta(b))$ is deterministic (as shown in part (b) of the proof below) which makes model B easier to analyze. 
We note here that, unlike the corresponding result in \cite{FS} for model B without disorder, we are unable to adapt to our disordered case the proof of Theorem 2 from \cite{bric} used in \cite{FS}. The proof in \cite{bric} relies on the weak convergence of $\mu_\L^{\rho_0}[\xi=0]$ to an infinite-volume gradient Gibbs measure $\mu[\xi=0]$ (which, due to the disorder, we were unable to show for $\mu^{\rho_0}_\L[\xi]$, but only for $\hat\mu_k^u[\xi]$, even for the periodic boundary conditions considered in \cite{bric}), and on the resulting Brascamp-Lieb inequality for the measure $\mu[\xi=0]$.

\begin{itemize}

\item [(a)]

We will first show the statement of the theorem for $u=0$, and then we will adapt the proof to the general $u\in\R^d$ case. For $u=0$, we will show that the random gradient Gibbs measure $\mu[\xi]$ constructed in Proposition \ref{existgibbs} satisfies (\ref{firstlimit}). For the general case $u\in\R^d$ we will follow the same approach as in \cite{FS} and use the fact that boundary conditions with definite
tilt $u$ are identical to boundary conditions $u=0$ for the shifted potential $V(\cdot+u_\alpha)$ for a bond in direction $e_\alpha$, where  $\alpha\in \{1,2,\ldots, d\}$. Thus an infinite-volume gradient Gibbs measure $\mu[\xi] $ with arbitrary expected tilt $u$ which satisfies Definition \ref{shiftcov1} is constructed from
the finite-volume gradient Gibbs measures with potential $V(\cdot+u_\alpha)$.

\vspace{2mm}

{\bf Step 1:}  Fix $\alpha\in \{1,2,\ldots, d\}$. We will show here that in order to prove (\ref{firstlimit}) for $u\in\R^d$, it is sufficient to prove that
\begin{equation}
\label{no1}
\liminf_{n\rightarrow\infty}\liminf_{k\rightarrow\infty}\frac{1}{k}\sum_{i=1}^k \frac{1}{|{\L}_{m_i}|}\sum_{w\in {\L}_{m_i}}\E\mu^{\rho_u}_{{\L}_{m_i}+w}[\xi]\left(\frac{1}{|\Lambda_n|}\sum_{x\in\Lambda_n}\eta(b_{x,\alpha})-u_\alpha\right)^2=0.
\end{equation}

\vspace{1.5mm}

We note first that since $\mu[\xi]$ satisfies the integrability assumption (\ref{intconda}), we have by a standard subadditivity argument (see, for example, \cite{Steele})
$$\lim_{|\Lambda|\rightarrow\infty}\left|\frac{1}{|\Lambda|}\sum_{x\in\Lambda}\eta(b_{x,\alpha})-u_\alpha\right|~~~~~\mbox{exists}~~\mu^u[\xi]-\mbox{a.s.}.$$
It follows that in order to show (\ref{firstlimit}), it suffices to show that for $\P$-a.s. $\xi$
\begin{equation}
\label{firstlimit1a}
\mu^u[\xi]\left(\lim_{n\rightarrow\infty}\left(\frac{1}{|\Lambda_n|}\sum_{x\in\Lambda_n}\eta(b_{x,\alpha})-u_\alpha\right)^2\right)=0.
\end{equation}
By Fatou's lemma, it follows that to show (\ref{firstlimit1a}) it is enough to prove that for $\P$-a.s. $\xi$
\begin{equation}
\label{firstlimit12a}
\liminf_{n\rightarrow\infty}{\mu}^{u}[\xi]\left(\frac{1}{|\Lambda_n|}\sum_{x\in\Lambda_n}\eta(b_{x,\alpha})-u_\alpha\right)^2=0,
\end{equation}
or equivalently 
\begin{equation}
\label{firstlimit120a}
\liminf_{n\rightarrow\infty}{\E\mu}^{u}[\xi]\left(\frac{1}{|\Lambda_n|}\sum_{x\in\Lambda_n}\eta(b_{x,\alpha})-u_\alpha\right)^2=0.
\end{equation}
By the lower semi-continuity of $\left(\frac{1}{|\Lambda_n|}\sum_{x\in\Lambda_n}\eta(b_{x,\alpha})\right)^2$ and by the weak convergence of $\hat{\mu}_k^u[\xi]$ to $\mu^u[\xi]$, we then have
\begin{eqnarray*}
\E\mu^u[\xi]\left(\frac{1}{|\Lambda_n|}\sum_{x\in\Lambda_n}\eta(b_{x,\alpha})-u_\alpha\right)^2&\le&\liminf_{k\rightarrow\infty}\E\hat{\mu}_k^u[\xi]\left(\frac{1}{|\Lambda_n|}\sum_{x\in\Lambda_n}\eta(b_{x,\alpha})-u_\alpha\right)^2\\
&=&\liminf_{k\rightarrow\infty}\frac{1}{k}\sum_{i=1}^k \frac{1}{|{\L}_{m_i}|}\sum_{w\in {\L}_{m_i}}\E\mu^{\rho_u}_{{\L}_{m_i}+w}[\xi]\left(\frac{1}{|\Lambda_n|}\sum_{x\in\Lambda_n}\eta(b_{x,\alpha})-u_\alpha\right)^2.
\end{eqnarray*}
Combining (\ref{firstlimit12a}) with the above, (\ref{no1}) follows.  

We will focus in Steps 2 and 3 below on estimating (\ref{no1}) in the particular case with $u=0$. Fix $m_i\in\N, x\in {\L}_{m_i}$ and $n\in\N$. 
We have 
\begin{eqnarray}
\label{BLsplit}
\lefteqn{\E\mu^{\rho_0}_{{\L}_{m_i}+w}[\xi]\left(\frac{1}{|\Lambda_n|}\sum_{x\in\Lambda_n}\eta(b_{x,\alpha})-u_\alpha\right)^2}\nonumber\\
&=&\E\left(\var_{\mu^{\rho_0}_{{\L}_{m_i}+w}[\xi]}\left(\frac{1}{|\Lambda_n|}\sum_{x\in\Lambda_n}\eta(b_{x,\alpha})-u_\alpha\right)\right)
+\Var\left({\mu^{\rho_0}_{{\L}_{m_i}+w}[\xi]}\bigg(\frac{1}{|\Lambda_n|}\sum_{x\in\Lambda_n}\eta(b_{x,\alpha})-u_\alpha\bigg)\right)\nonumber\\
&&+\left(\E\mu^{\rho_0}_{{\L}_{m_i}+w}[\xi]\left(\frac{1}{|\Lambda_n|}\sum_{x\in\Lambda_n}\eta(b_{x,\alpha})-u_\alpha\right)\right)^2.
\end{eqnarray}
We will estimate in Steps 2 and 3 below each of these three terms above separately for the $u=0$ case.

\vspace{2mm}

{\bf Step 2:} We will prove in this step that for all $m_i\in\N,x, w\in\zd$, we have
\begin{equation}
\label{lala}
\E\nu^{0}_{{\L}_{m_i}+w\setminus \{w\}}[\xi]\left(\phi(x)\right)=0,
\end{equation}
where we denoted by $\nu^{0}_{ {\L}_{m_i}+w\setminus \{0\}}[\xi]$ the Gibbs measure with $0$ boundary conditions \textit{outside $\L_{m_i+w}$} and \textit{at $w$}. Since by (\ref{19})
$$\E\mu^{\rho_0}_{{\L}_{m_i}+w}[\xi]\left(\sum_{x\in\Lambda_n}\eta(b_{x,\alpha})\right)=\sum_{x\in\Lambda_n}\E\nu^{0}_{{\L}_{m_i}+w\setminus\{w\}}[\xi]\left(\phi(x+e_\alpha)-\phi(x)\right),$$
this will imply that the third term on the right-hand side in (\ref{BLsplit}) is equal to $0$.

\vspace{1.5mm}

To show (\ref{lala}) we will take advantage of the symmetry of $V$. More precisely, by means of the change of variables $\phi(y)\rightarrow-\phi(y)$, $y\in {\L}_{m_i}+w$, we have
$$\nu^{0}_{{\L}_{m_i}+w\setminus \{w\}}[\xi](\phi(x))=-\nu^{0}_{{\L}_{m_i}+w\setminus \{w\}}[-\xi](\phi(x)).$$
Using now the independence of the disordered random fields $(\xi(x))_{x\in\Z^d}$ and the symmetry of their distribution, we get in the above
$$\E\nu^{0}_{{\L}_{m_i}+w\setminus \{w\}}[\xi](\phi(x))=-\E\nu^{0}_{{\L}_{m_i}+w\setminus \{w\}}[-\xi](\phi(x))=-\E\nu^{0}_{{\L}_{m_i}+w\setminus \{w\}}[\xi](\phi(x)),$$
from which (\ref{lala}) immediately follows.

\vspace{2mm}

{\bf Step 3:} We will estimate here the first two terms in (\ref{BLsplit}).

\vspace{1.5mm}

We need only consider the case with $\Lambda_n\cap \Lambda_{m_i+x}\neq \emptyset$ as otherwise (\ref{BLsplit}) is $0$ due to the boundary conditions. By the Brascamp-Lieb inequality (\ref{bl11}), we have for the first term on the right-hand side in (\ref{BLsplit})
\begin{equation}
\label{BLsplit1}
\var_{\mu^{\rho_0}_{{\L}_{m_i}+w}[\xi]}\left(\frac{1}{|\Lambda_n|}\sum_{x\in\Lambda_n}\eta(b_{x,\alpha})\right)
\le \frac{1}{C_1}\mu^{\rho_0}_{G,{\L}_{m_i}+w}[\xi=0]\left(\frac{1}{|\Lambda_n|}\sum_{x\in\Lambda_n}\eta(b_{x,\alpha})\right)^2.
\end{equation}
In order to estimate this further, we will need to introduce first some notation.
 Let $\Lambda_{m_i+w,n}:=\L_{m_i+w}\cap\L_n$, let $\partial\Lambda^{+}_{m_i+w,n}$ be the boundary of $\Lambda_{m_i+w,n}$ and let $\partial\Lambda^{-}_{m_i+w,n}:=\{a\in \Lambda_{m_i+w,n}~|~\exists y\in \partial\Lambda^{+}_{m_i+w,n}~~\mbox{such that}~~|a-y|=1\}$. We note here that $\left|\partial\Lambda^{-}_{m_i+w,n}\right|\le (2n)^{d-1}$, which fact will be used a few times in the proof. Taking account of boundary conditions, of term cancellations and of Proposition \ref{properties} (ii), we have for the right-hand side of (\ref{BLsplit1}) 
 \begin{eqnarray}
 \label{blsplit2}
 \mu^{\rho_0}_{G,{\L}_{m_i}+w}[\x=0]\left(\frac{1}{|\Lambda_n|}\sum_{x\in\Lambda_n}\eta(b_{x,\alpha})\right)^2&\le& \nu^{0}_{G, {\L}_{m_i}+w\setminus \{w\}}[\xi=0]\left(\frac{1}{|\Lambda_n|}\sum_{y\in \partial\Lambda^{-}_{m_i+x,n}}\phi(y)\right)^2\nonumber\\
 &\le&\frac{1}{(2n)^{d+1}}\sum_{y\in \partial\Lambda^{-}_{m_i+w,n}} \nu^{0}_{G, {\L}_{m_i}+w\setminus \{w\}}[\xi=0]\left(\phi(y)\right)^2\nonumber\\
 &\le&\frac{1}{(2n)^{d+1}}\sum_{y\in \partial\Lambda^{-}_{m_i+w,n}} G_{{\L}_{m_i}+w}(y,y)\le\frac{C(d)}{n^2},
 \end{eqnarray}
 for some constant $C(d)>0$, independent of $m_i, n, \xi, w$ and $x$, and where $ \nu^{0}_{G, {\L}_{m_i}+w\setminus \{w\}}[\xi=0]$ is a Gaussian Gibbs measure with $0$ boundary conditions outside $\L_{m_i+w}$ and \textit{at $w$}. We note here that the pinning of the measure at $w$ plays no role for model A in the computations above, but will be \textit{crucial} in the corresponding computations for bounding the variance in (\ref{blsplit2}) for model B in $d=1,2$. We will next estimate the second term on the right-hand side of (\ref{BLsplit}). 
By means of Proposition \ref{golemma} and by using the fact that $(\xi(x))_{x\in\zd}$ are i.i.d., we have 
\begin{eqnarray}
\label{goeqn}
\lefteqn{\Var\left({\mu^{\rho_0}_{{\L}_{m_i}+w}[\xi]}\bigg(\frac{1}{|\Lambda_n|}\sum_{x\in\Lambda_n}\eta(b_{x,\alpha})\bigg)\right)}\nonumber\\
&\le&\Var(\xi(0))\sum_{z\in{\L}_{m_i}+w}\E\left(\sup_{\xi(z)}\cov^2_{\nu^{0}_{{\L}_{m_i}+w\setminus \{w\}}[\xi]}\left(\phi(z),\frac{1}{|\Lambda_n|}\sum_{x\in\Lambda_n}\eta(b_{x,\alpha})\right)\right).
\end{eqnarray}
To bound (\ref{goeqn}) we will consider separately the cases $d\ge 5$ and the critical cases $d=3,4$.

\vspace{2mm}

\textbf{(i) Case $d\ge 5$}. Then we have from (\ref{goeqn}) and (\ref{19})
\begin{eqnarray}
\label{gogo}
\lefteqn{\Var\left({\mu^{\rho_0}_{{\L}_{m_i}+w}[\xi]}\bigg(\frac{1}{|\Lambda_n|}\sum_{x\in\Lambda_n}\eta(b_{x,\alpha})\bigg)\right)}\nonumber\\
&\le&\Var(\xi(0))\sum_{z\in{\L}_{m_i}+w}\E\bigg(\sup_{\xi(z)}\bigg( \frac{1}{|\Lambda_n|}\sum_{y\in \partial\Lambda^{-}_{m_i+w,n}}\cov_{\nu^{0}_{{\L}_{m_i}+w\setminus \{w\}}[\xi]}\bigg(\phi(z), \phi(y)\bigg)\bigg)^2\bigg)\nonumber\\
&\le&\frac{\Var(\xi(0))}{n^{d+1}}\sum_{y\in \partial\Lambda^{-}_{m_i+w,n}}\sum_{z\in{\L}_{m_i}+w}\E\bigg(\sup_{\xi(z)}\cov^2_{\nu^{0}_{{\L}_{m_i}+w\setminus \{w\}}[\xi]}\bigg(\phi(z),  \phi(y)\bigg)\bigg)\nonumber\\
&\le&\frac{\Var(\xi(0))}{n^{d+1}}\sum_{y\in \partial\Lambda^{-}_{m_i+w,n}}\sum_{z\in{\L}_{m_i}+w}\frac{C'(d)}{]|y-z|[^{2d-4}}\nonumber\\
&\le&\frac{\Var(\xi(0))}{n^{d+1}}\sum_{y\in \partial\Lambda^{-}_{m_i+w,n}}C''(d)=\frac{\Var(\xi(0))C''(d)}{n^{2}},
\end{eqnarray}
where for the second inequality we used $(\sum_{i\in I} a_i)^2\le |I| \sum_{i\in I} a_i^2,$ which trivially holds for any finite set $I\subset \zd$ and for any $(a_i)_{i\in I}\in\R^I$, and for the third inequality we used the random walk representation estimates from Proposition \ref{propg} (ii). Note that by Proposition \ref{propg} (ii), $C'(d), C''(d)>0$ are independent of $m_i,x,n,w$ and of the disorder $\xi$. Combining (\ref{gogo}) with (\ref{no1}), (\ref{BLsplit}) and (\ref{lala}) proves the theorem in this case.

\vspace{2mm}

\textbf{ (ii) Case $d=3,4$}. In this case, estimating the sum on the right-hand side of (\ref{goeqn}) by the suboptimal estimates in (\ref{gogo}) would lead to a bound depending on $m_i$ if $|\Lambda_n|$ and $|\Lambda_{m_i+x}|$ are not of the same order. Since we need to look at estimates for all boxes, due to the fact that we average over them in (\ref{no1}), we will proceed as follows. For  $\Lambda_{m_i+x}\subset \Lambda_{2n}$ we will estimate the variance as in (\ref{gogo}) and we have
\begin{eqnarray}
\label{gogo11}
\mathbb{\Var}\left({\mu^{\rho_0}_{{\L}_{m_i}+w}[\xi]}\bigg(\frac{1}{|\Lambda_n|}\sum_{x\in\Lambda_n}\eta(b_{x,\alpha})\bigg)\right)&\le&\frac{\Var(\xi(0))}{n^{d+1}}\sum_{y\in \partial\Lambda^{-}_{m_i+w,n}}\sum_{z\in{\L}_{2n}}\frac{C'(d)}{]|y-z|[^{2d-4}}\nonumber\\
&\le&\frac{n\Var(\xi(0))}{n^{d+1}}\sum_{y\in \partial\Lambda^{-}_{m_i+w,n}}C'''(d)\nonumber\\
&=&\frac{\var(\xi(0))C'''(d)}{n},
\end{eqnarray}
where $C'(d), C'''(d)>0$ are independent of $m_i,x,n$ and of the disorder $\xi$.
For $\Lambda_{2n}\subset \Lambda_{m_i+w}$ we have
\begin{multline}
\label{gogo1}
\Var\left({\mu^{\rho_0}_{{\L}_{m_i}+w}[\xi]}\bigg(\frac{1}{|\Lambda_n|}\sum_{x\in\Lambda_n}\eta(b_{x,\alpha})\bigg)\right)\\\le\Var(\xi(0))\sum_{z\in{\L}_{2n}}\E\bigg(\sup_{\xi(z)}\cov^2_{\nu^{0}_{{\L}_{m_i}+w\setminus \{w\}}[\xi]}\bigg(\phi(z), \frac{1}{|\Lambda_n|}\sum_{y\in \partial\Lambda^{-}_{m_i+x,n}} \phi(y)\bigg)\bigg)\\
+\Var(\xi(0))\sum_{z\in {\L}_{m_i+w}\setminus {\L}_{2n}}\E\bigg(\sup_{\xi(z)}\cov^2_{\nu^{0}_{{\L}_{m_i}+w\setminus \{w\}}[\xi]}\bigg(\phi(z),\frac{1}{|\Lambda_n|}\sum_{x\in\Lambda_n}\eta(b_{x,\alpha}) \bigg)\bigg).
\end{multline}
The first term on the right-hand side above can be estimated as in (\ref{gogo11}); recalling (\ref{dgirep1a}), we have for the second term 
\begin{eqnarray}
\label{gogo2}
\lefteqn{\sum_{z\in {\L}_{m_i+w}\setminus {\L}_{2n}}\E\bigg(\sup_{\xi(z)}\cov^2_{\nu^{0}_{{\L}_{m_i}+w\setminus \{w\}}[\xi]}\bigg(\phi(z),\frac{1}{|\Lambda_n|}\sum_{x\in\Lambda_n}\eta(b_{x,\alpha}) \bigg)\bigg)}\nonumber\\
&\le& \frac{1}{|\Lambda_n|}\sum_{x\in\Lambda_n}\sum_{z\in {\L}_{m_i+w}\setminus {\L}_{2n}}\E\bigg(\sup_{\xi(z)}\cov^2_{\nu^{0}_{{\L}_{m_i}+w\setminus \{w\}}[\xi]}\bigg(\phi(z),\eta(b_{x,\alpha}) \bigg)\bigg)\nonumber\\
&=&\frac{1}{|\Lambda_n|}\sum_{x\in\Lambda_n}\sum_{ z\in {\L}_{m_i+w}\setminus {\L}_{2n}}\E\bigg(\sup_{\xi(z)}\left(\nu^{0}_{{\L}_{m_i}+w\setminus \{w\}}[\xi]\left(\int_0^\infty\nabla_\alpha p^{\nabla\phi}_{{\L}_{m_i}+w}(0,x,t, z)\rmd t\right)\right)^2\bigg)\nonumber\\
&=&\frac{1}{|\Lambda_n|}\sum_{x\in\Lambda_n}\sum_{z\in {\L}_{m_i+w}\setminus {\L}_{2n}}\E\bigg(\sup_{\xi(z)}\left(\nu^{0}_{{\L}_{m_i}+w\setminus \{w\}}[\xi]\left( \nabla_\alpha g^{\nabla\phi}_{{\L}_{m_i}+w}(x, z)\right)\right)^2\bigg),
\end{eqnarray}
where for the first equality we used Proposition \ref{dgirep}, and where $\nabla_\alpha p^{\nabla\phi}_{{\L}_{m_i}+w}(0,x,t, z):=p^{\nabla\phi}_{{\L}_{m_i}+w}(0,x,t, z)-p^{\nabla\phi}_{{\L}_{m_i}+w}(0,x+e_\alpha,t, z)$, with a similar definition for $\nabla_\alpha g^{\nabla\phi}_{{\L}_{m_i}+w}(x,z)$. Note now that for all $z\in  {\L}_{m_i+w}\setminus {\L}_{2n}$ and $x\in {\L}_n$ we have $|x-z|\ge n$. 

For $d=4$, it follows now easily from  Proposition \ref{propg} (iv) that the quantity in (\ref{gogo2}) is bounded by $C(4)/n^\delta$, for some $C(4)$ which is independent of $m_i,x,w$ and $n$. Combining (\ref{no1}), (\ref{BLsplit}),  (\ref{gogo11}), (\ref{gogo1}), (\ref{gogo2}), (\ref{gogo4}) and (\ref{lala}) proves the theorem for $d=4$.

We focus next on the more delicate $d=3$ case. Since the estimates from Proposition \ref{propg} (ii) and (iv) are too weak for $d=3$ to give us a bound in (\ref{gogo2}) which is independent of $m_i$, we will re-write (\ref{gogo2}) in a form in which we can use (\ref{gogo30a}). As a result, we need to work under the more restrictive assumption (\ref{gap}) on the disorder, which allows us to get rid of the supremum in  (\ref{gogo2}). Note first that 
$${\L}_{m_i}+w\setminus {\L}_{2n}\subset\cup_{j=1}^{1+\left[\log(\frac{{3m_i}}{n})\right]} \left({\L}_{2^{j+1}n}\setminus {\L}_{2^{j}n}\right),$$ 
with $[x]$ the integer part of $x$. In particular, for all $ z\in  {\L}_{2^{j+1}n}\setminus {\L}_{2^{j}n}$ and $x\in {\L}_n$, $j\ge 1$, we have $|x-z|\ge 2^{j-1} n$. We have now in view of (\ref{gogo2}), (\ref{no1}) and of $g^{\nabla\phi}_{{\L}_{m_i}+w}(x,z)= g^{\tau_{-z}(\nabla\phi)}_{{\L}_{m_i}+w-z}(x-z,0)$ (which follows from (\ref{dgirep1a}) by the shift $\phi(v)\rightarrow\phi(v-z),v\in\zd$)
\begin{eqnarray}
\label{gogo4}
\lefteqn{\frac{1}{|\L_{m_i}|}\sum_{w\in\L_{m_i}}\sum_{z\in {\L}_{m_i+w}\setminus {\L}_{2n}}\E\bigg(\cov^2_{\nu^{0}_{{\L}_{m_i}+w\setminus \{w\}}[\xi]}\bigg(\phi(z),\eta(b_{x,\alpha}) \bigg)\bigg)}\nonumber\\
&=&\frac{1}{|\L_{m_i}|}\sum_{w\in\L_{m_i}}\sum_{z\in {\L}_{m_i+w}\setminus {\L}_{2n}} \E\bigg(\left(\nu^{0}_{{\L}_{m_i}+w-z\setminus \{w-z\}}[\tau_{-z}\xi]\left( \nabla_\alpha g^{\nabla\phi}_{{\L}_{m_i}+w-z}(x-z, 0)\right)\right)^2\bigg)\nonumber\\
&=& \frac{1}{|\L_{m_i}|}\sum_{w\in\L_{m_i}}\sum_{j=1}^{1+\left[\log(\frac{{3m_i}}{n})\right]}\sum_{z\in {\L}_{2^{j+1}n}\setminus {\L}_{2^jn}}\E\bigg(\left(\nu^{0}_{{\L}_{m_i}+w-z\setminus \{w-z\}}[\tau_{-z}\xi]\left( \nabla_\alpha g^{\nabla\phi}_{{\L}_{m_i}+w-z}(x-z, 0)\right)\right)^2\bigg)\nonumber\\
&\le& \frac{1}{|\L_{m_i}|}\sum_{v\in\L_{2m_i}}\sum_{j=1}^{1+\left[\log(\frac{{3m_i}}{n})\right]}\sum_{w,z\in\L_{m_i}: w-z=v\atop z\in {\L}_{2^{j+1}n}\setminus {\L}_{2^jn}}\E\bigg(\nu^{0}_{{\L}_{m_i}+w-z\setminus \{w-z\}}[\xi]\left( \nabla_\alpha g^{\nabla\phi}_{{\L}_{m_i}+w-z}(x-z, 0)\right)^2\bigg)\nonumber\\
&\le&\frac{\tilde{C}}{|\Lambda_{m_i}|}\sum_{v\in\Lambda_{2m_i}}\sum_{j=1}^{1+\left[\log(\frac{{3m_i}}{n})\right]} \frac{1}{2^{j-1}n}\le\frac{C'}{n},
\end{eqnarray}
for some $C'>0$ independent of $m_i,x,w$ and $n$, and where for the first inequality we used the fact that $(\xi(y))_{y\in\zd}$ are i.i.d., and for the second inequality we used (\ref{gogo30a}) from Proposition \ref{propg}. Combining now (\ref{no1}), (\ref{BLsplit}), (\ref{gogo11}), (\ref{gogo1}), (\ref{gogo2}), (\ref{gogo4}) and (\ref{lala}) proves the theorem.

\vspace{2mm}

{\bf Step 4:} We will show here (\ref{firstlimit}) for the general $u\in\R^d$ case.

With the usual notations, let us define the shifted measure
$$\nu_{\mbox{\shift}, \Lambda}^\psi[\xi ](\ormd\phi):=\frac{1}{{Z}_{\mbox{\shift},\Lambda}^{\psi}[\xi]}e^{- \frac{1}{2}\sum_{x\in\Lambda, y\in\Lambda\cup\partial\Lambda\atop |x-y|=1}V(\phi(x)-\phi(y)-\langle u,x-y \rangle)+ \sum_{x\in\Lambda}\xi(x)\phi(x)}\rmd\phi_\Lambda\delta_\psi(\ormd\phi_{{\Z}^d\setminus\Lambda}),$$
and let $\mu_{\mbox{\shift}, \Lambda}^\rho[\xi](\ormd\eta)$ be the corresponding finite-volume gradient Gibbs measure on $\chi$ such that Definition \ref{finvolgrad} is satisfied. Let
$$\hat{\mu}^u_{{\shift,k}[\xi]}:=\frac{1}{k}\sum_{i=1}^k {\bar\mu}^{u}_{\shift, \L_{m_{i}}}[\xi],$$
where $\bar{\mu}^u_{\shift, \L_{m_{i}}}$ is defined as in (\ref{heldfixed}). 
We can now reason as in \cite{CK} to show that $\hat{\mu}^u_{{\shift,k}[\xi]}$ converges weakly to a shift-covariant gradient Gibbs measure ${\mu}^u_{\shift}[\xi]$ which satisfies Definition \ref{shiftcov1}. That is, we will first show as in Proposition 3.6 from \cite{CK} that $$\P^u_{{\shift},\L}(\ormd\phi)
:=\left(\int\P(\ormd\xi)\bar\mu^{u}_{{\shift},\L}[\xi]\right)(\ormd\phi)$$
satisfies for some $K>0$, uniformly in $x_0,y_0\in\zd$, the estimate 
\begin{equation}
\label{tight}
\limsup_{N\uparrow \infty} \P^u_{{\shift},\L_N}\left[(\phi(x_0)-\phi(y_0)-u\cdot(x_0 -y_0))^2\right] \leq K.
\end{equation}
The key idea is to perform in (\ref{tight}) the change of variables $\phi(x)\rightarrow\tilde\phi(x)+x\cdot u, x\in\zd$, which shifts $\P^u_{{\shift},\L}\left[(\phi(x_0)-\phi(y_0)-u\cdot(x_0 -y_0))^2\right]$ to $\P^0_{\L}\left[(\phi(x_0)-\phi(y_0))^2\right] 
:=\int\P(\ormd\xi)\bar\mu^{0}_{\L}[\xi](\phi(x_0)-\phi(y_0))^2$.
By (\ref{tight}) the sequence of measures $\P^u_{{\shift},\L_N}$ is tight. By the same arguments as in Proposition 3.8 from \cite{CK} we can show that $\hat{\mu}^u_{{\shift,k}}[\xi]$ converges weakly to a shift-covariance gradient Gibbs measure $\tilde{\mu}^u_{\shift}[\xi]$ satisfying Definition \ref{shiftcov1}. Moreover, $\tilde{\mu}^u_{\shift}[\xi]$  can be shown as in Step 2 above, by the same change of variables $\phi(x)\rightarrow\tilde\phi(x)+x\cdot u, x\in\zd$, to have expected tilt $u$. 

The proof of (\ref{firstlimit}) now follows the same reasoning as in Steps 1,2 and 3 above. 

\item [(b)] For $u=0$ we have by symmetry of $V_{(x,y)}$ that for all $m_i\in\N,x, w\in\zd$,
$\nu^{0}_{{\L}_{m_i}+w\setminus\{w\}}[\omega]\left(\phi(x)\right)=0$.
Therefore, the proof reduces to finding an upper bound for
$$\var_{\nu^{0}_{{\L}_{m_i}+w\setminus\{w\}}[\omega]}\left(\frac{1}{|\Lambda_n|}\sum_{x\in\Lambda_n}\eta(b_{x,\alpha})\right)^2,$$
which can be easily done by the Brascamp-Lieb inequality (\ref{bl11}) and (for the critical cases $d=1,2$) also by the estimates from (\ref{varloc}). The extension to $u\in\R^d$ follows as in Step 4 above.
\end{itemize}
$\Cox$

\begin{rem}
\label{rema}
\begin{enumerate}
\item [(a)] Note that (\ref{firstlimit}) (respectively (\ref{firstlimit11})) implies that $\mu[\x]$ (respectively $\mu[\omega]$) has expected tilt $u$, that is
\begin{equation*}
 \E\left(\int\mu^u[\xi](\ormd\eta)\eta(b)\right)=\langle u,y_b-x_b \rangle ~\mbox{for all bonds}~b=(x_b,y_b)\in(\Z^d)^*.
\end{equation*}
\item [(b)] Property (\ref{firstlimit}) (respectively property (\ref{firstlimit11})) is not preserved under a convex combination of measures with different expected tilts. That is, let $u_1\in\R^d$, $u_2\in\R^d$ and $a\in [0,1]$. Let $\mu^{u_1}[\xi]$ and $\mu^{u_2}[\xi]$  be two measures defined as in Definition \ref{shiftcov1}, with expected tilts $u_1$ and $u_2$, which
satisfy (\ref{firstlimit}) for $\P$-almost every $\xi$. Then $a\mu^{u_1}[\xi] +(1-a)\mu^{u_2}[\xi]$ need not satisfy (\ref{firstlimit}),
even though $\E(a\mu^{u_1}[\xi](\eta(b)) +(1-a)\mu^{u_2}[\xi](\eta(b)))= \langle au_1+(1-a)u_2, y_b-x_b \rangle ~\mbox{for all bonds}~b=(x_b,y_b)\in(\Z^d)^*$. 
\item [(c)]  For model B, our proof can be applied to a class of non-convex potentials at \textbf{all temperatures}, since for (\ref{firstlimit11}) to hold, we only need an upper bound on the variance, uniform in the size of the box. This can be done by an extension of the Brascamp-Lieb inequality to a class of non-convex potentials, as shown for example in Proposition A.2 from \cite{H}. For potentials \textit{without} disorder, in view of the ergodic decomposition of shift-invariant Gibbs measures (see, for example, Chapter 14 from \cite{giorgii} for more on this), (\ref{firstlimit11})  implies existence of \textbf{ergodic, extremal} gradient Gibbs measures with given tilt for a certain class of non-convex potentials at all temperatures, which class includes the potential studied in \cite{BK}. 
\end{enumerate}
\end{rem}

\section{Dynamical method: coupling gradient Gibbs measures with given averaged tilt for the same disorder and same dynamics}

The main result proved in this section is Theorem \ref{uniq}. The proof will be done in two steps. First, in subsection 4.1 we will prove in Theorem \ref{extremeuniq} a statement of uniqueness of shift-covariant gradient Gibbs measure with direction-averaged tilt. The proof of Theorem \ref{extremeuniq} relies on a far from trivial adaptation of the method of Funaki and Spohn in Theorem 2.1 from \cite{FS}, to obtain uniqueness of the gradient Gibbs measure averaged over the disorder with direction-averaged tilt. Proposition \ref{kom} allows us to transform this into a statement of uniqueness of the corresponding quenched gradient Gibbs measure with direction-averaged expected tilt.  Then we will upgrade this statement to the one in Theorem \ref{uniq} by using the quenched uniqueness result in Theorem \ref{extremeuniq} and a proof by contradiction argument.



\subsection{Uniqueness of gradient Gibbs measure with given direction-averaged tilt}

Before we state the main result of this section, Theorem \ref{extremeuniq} below, we will introduce the dynamics which govern the $\phi$- and the $\eta$-fields. Because of long-range dependence, Dobrushin type methods do not seem to work for the uniqueness problem for gradient models with or without disorder, which is why both in \cite{FS} and in our proof the dynamics is used to help establish the result. We assume that the dynamics of the height variables $\phi_t=\{\phi_t(y)\}_{y\in\zd}$ are generated by the following family of SDEs:
\begin{itemize}
\item [(A)] For model (A), we have for all $\xi\in\Omega$
\begin{equation}
\label{sde}
\rmd\phi_t(y)=-\sum_{x\in\zd,||x-y||=1}V'(\phi_t(x)-\phi_t(y))\rmd t+\xi(y)\rmd t+ \sqrt{2}d W_t(y),~~y\in\zd,
\end{equation}
where $\{W_t(y),y\in\zd\}$ is a family of independent Brownian motions.  
The dynamics for the height differences $\eta_t=\{\eta_t(b)\}_{b\in (\zd)^*}$ are then determined for all $b\in(\zd)^*$ by
\begin{equation}
\label{sde1}
\rmd\eta_t(b)=-\sum_{b'\in (\zd)^*: x_{b'}=x_b}V'(\eta(b'))\rmd t+\xi(x_b)\rmd t+ \sqrt{2}d W_t(b),~~b\in(\zd)^*,
\end{equation}
where $W_t(b):=W_t(x_b)-W_t(y_b)$.
\item [(B)] For model (B), we have for all $\omega\in\Omega$
\begin{equation}
\label{sdeb}
\rmd\phi_t(y)=-\sum_{x\in\zd,||x-y||=1}(V^{\omega}_{\langle x,y\rangle})'(\phi_t(x)-\phi_t(y))\rmd t+ \sqrt{2}d W_t(y),~~y\in\zd,
\end{equation}
where $\{W_t(y),y\in\zd\}$ is a family of independent Brownian motions. 
The dynamics for the height differences $\eta_t=\{\eta_t(b)\}_{b\in (\zd)^*}$ are then determined by
\begin{equation}
\label{sdeb1}
\rmd\eta_t(b)=-\sum_{b'\in (\zd)^*: x_{b'}=x_b}(V_{b'}^{\omega})'(\eta(b'))\rmd t+ \sqrt{2}d W_t(b),~~b\in(\zd)^*.
\end{equation}
\end{itemize}
Due to the conditions on the potentials in both models (A) and (B) and to the second moments assumption on the disorder in model (A),  there is global Lipschitz continuity in $\chi_r, r>0,$ on the drift part of the SDEs. Then, as a consequence of an infinite version of the Yamada-Watanabe result of existence and uniqueness of strong solutions to SDEs (as stated, for example, in \cite{GM}), one can show that (\ref{sde1}) and (\ref{sdeb1}) 
have a unique $\chi_r$-valued continuous strong solution starting at $\eta_0=\eta\in\chi$. 

\vspace{2mm}

Let ${\cal P}(\chi)$ be the set of all probability
measures on $\chi$ and let ${\cal P}_2(\chi)$ be those $\mu\in {\cal P}(\chi)$ satisfying $E_{\mu}[|\eta(b)|^2]<\infty$ for each $b\in (\zd)^*$. For $r>0$, recall the definition of $\chi_{r}$ as given in Subsection \ref{nablagibbs}. The set ${\cal P}(\chi_r), r>0$, is defined correspondingly and ${\cal P}_2(\chi_r)$ stands for the set of all $\mu\in {\cal P}(\chi_r)$ such that $E_{\mu}[|\eta|^2_r]<\infty$.

\vspace{3mm}

We are now ready to state the main result of this section: 
\begin{thm} 
\label{extremeuniq}  
Let $u\in\R^d$. Recall that for all $\alpha\in\{1,2,\ldots,d\}$ we defined
$$E_\alpha:=\{ \eta~|~\lim_{|\Lambda|\rightarrow\infty}\frac{1}{|\Lambda|}\sum_{x\in\Lambda}\eta(b_{x,\alpha})=u_\alpha\},$$
along the sequence of volumes with $b_{x,\alpha}:=(x+ e_\alpha,x)\in (\zd)^*$. 
\begin{enumerate}
\item [(a)] \textbf{(Model A)} Let $d\ge 3$. Assume that $V$ satisfies (\ref{tag22}) and that $(\xi(x))_{x\in\zd}$ have symmetric distribution. For $d=3$ we will also assume that the distribution of $\xi(0)$ satisfies (\ref{gap}). Then there exists \textbf{at most one} $\P$-almost surely shift-covariant measure $\xi\rightarrow\mu[\xi]$, $\mu[\xi] \in {\cal P}(\chi)$, stationary for the SDE (\ref{sde1}), which satisfies for $\P$-almost every $\xi$
\begin{equation*}
\label{firstlimit12}
\mu^u[\xi](E_\alpha)=1,~\alpha\in\{1,2,\ldots,d\},
\end{equation*} 
and which satisfies the integrability condition
\begin{equation*}
\label{intcond1}
\E\int\mu^u[\xi](\ormd\eta)(\eta(b))^2<\infty~\mbox{for all bonds}~b\in (\zd)^*.
\end{equation*}

\item[(b)] \textbf{(Model B)} Let $d\ge 1$. Assume that for $\P$-almost every $\omega$, $V^\omega_{(x,y)}$ satisfies (\ref{tag23}) uniformly in the bonds $(x,y)$. Then 
there exists \textbf{at most one} $\P$-almost surely shift-covariant measure $\omega\rightarrow\mu[\omega]$, $\mu[\omega] \in {\cal P}(\chi)$, stationary for the SDE (\ref{sdeb1}), which satisfies for $\P$-almost every $\omega$
\begin{equation*}
\label{secondlimit12}
\mu^u[\omega](E_\alpha)=1,~\alpha\in\{1,2,\ldots,d\},
\end{equation*} 
and which satisfies the integrability condition
\begin{equation*}
\label{intcond2}
\E\int\mu^u[\omega](\ormd\eta)(\eta(b))^2<\infty~\mbox{for all bonds}~b\in (\zd)^*.
\end{equation*}
\end{enumerate}
\end{thm}
We will only do the proof of Theorem \ref{extremeuniq} for model (A), as the proof for model (B) follows similarly. We will prove Theorem \ref{extremeuniq} by coupling
techniques. We will follow the same line of argument as in \cite{FS}, by introducing dynamics on the
gradient field. However as we already emphasized, we do not have shift-invariance and ergodicity of the quenched measure as there is for the measure without disorder in \cite{FS}, which complicates matters considerably in our case.

The basic idea is as follows. Take two random gradient Gibbs measures (potentially different) with the same expected tilt; we know they are both invariant under the same stochastic dynamics.   
Take two initial realizations of field configurations corresponding to these gradient measures,  
and compute the change of distance between the evolved configurations of fields 
between time $0$ and a time $T$ as an integral over a time-derivative.
This time-derivative can be related to the distance of time-evolved gradient configurations
corresponding to the two initial conditions by means of the uniform strict convexity of the potential.
Taking expectations 
over the initial configurations and over the coupling dynamics, and then dividing the equation 
by large $T$ so that the contributions from time zero and $T$ drop out, 
one produces a coupling between the two shift-covariant gradient Gibbs measures. The expectation w.r.t. a certain averaged version of
this coupling measure 
becomes  arbitrarily small when $T$ is large.  
This proves the desired equality of the gradient Gibbs measures.


Formally, the proof of Theorem \ref{extremeuniq} is based on a coupling lemma, Lemma \ref{2.1} below; a key ingredient for the coupling lemma is a bound
on the distance between two measures evolving under the same dynamics. The main ingredients needed to prove the lemma are Theorem \ref{directilt}, a non-standard ergodic theorem for the measure averaged over the disorder (see (\ref{ergtheo}) below), the proof of uniqueness of the Gibbs measure averaged over the disorder from Lemma \ref{annealeduniq}, exploiting the rapid decay of the norm $\|\eta\|_r, r>0$, and Proposition \ref{kom} below (for a proof see Proposition 1a from \cite{KOM}). 
\begin {prop}
\label{kom}
If $(\zeta_n)_{n\in\N}$ is a sequence of real-valued random variables with $\lim\inf_{n\rightarrow\infty}\E(|\zeta_n|)<\infty$, there exists a subsequence $\{\theta_n\}_{n\in\N}$ of the sequence $\{\zeta_n\}_{n\in\N}$ and an integrable random variable $\theta$ such that for any arbitrary subsequence $\{\tilde{\theta}_n\}_{n\in\N}$ of the sequence $\{\theta_n\}$, we have almost surely that
$$\lim_{n\rightarrow\infty}\frac{\tilde{\theta}_1+\tilde{\theta}_2+\ldots+\tilde{\theta}_n}{n}=\theta.$$
\end{prop}



\noindent \textbf{Coupling Argument}

Take $u\in\R^d$. Suppose that there exist two shift-covariant measures $\xi\rightarrow\mu[\xi], \xi\rightarrow\bar\mu[\xi]$, $\mu[\xi],\bar\mu[\xi] \in {\cal P}(\chi)$, stationary for the SDE (\ref{sde1}), which satisfy for $\P$-almost every $\xi$
\begin{equation*}
\mu[\xi](E_\alpha)=1,~\bar\mu[\xi](E_\alpha)=1,~~\alpha\in\{1,2,\ldots,d\},
\end{equation*} 
and which satisfy the integrability condition
\begin{equation*}
\E\int\mu[\xi](\ormd\eta)(\eta(b))^2<\infty, \E\int\bar\mu[\xi](\ormd\eta)(\eta(b))^2<\infty,~~\mbox{for all bonds}~b\in (\zd)^*.
\end{equation*}
Note that $\E\int\mu[\xi](\ormd\eta), \E\int\bar\mu[\xi](\ormd\eta)$ are supported on ${\cal P}_2(\chi_r)$, for every $r>0$. We also note that one can show by means of Kolmogorov's characterization of reversible diffusions (see, for example, Corollary 1 in \cite{RRS} for the statement)  that every shift-covariant gradient Gibbs measure  $\xi\rightarrow \mu[\xi]$, defined as in Definition \ref{shiftcov1}, is reversible for the SDE (\ref{sde1}). (For the definition and proof of reversibility of Gibbs measures, see Proposition 3.1 in \cite{FS}.)  Moreover, the existence of such a shift-covariant gradient Gibbs measure satisfying the remaining conditions in Theorem \ref{extremeuniq} (a) is assured by Theorem \ref{directilt}(a).

\vspace{2mm}

For each fixed $\xi\in\Omega$, we construct two independent $\chi_{r}$-valued 
random variables $\eta=\{\eta(b)\}_{b\in (\zd)^*}$ and ${\bar{\eta}}=\{{\bar{\eta}}(b)\}_{b\in (\zd)^*}$ on a common probability space $(\Upsilon,{\cal L},\mathbb{Q}[\xi])$ in such a manner that $\eta$ and 
${\bar{\eta}}$ are distributed by $\mu[\xi]$ and ${\bar{\mu}}[\xi]$ under $\mathbb{Q}[\xi]$, respectively. We define $\phi_0=\phi^{\eta,0}$ and $\bar{\phi}_0=\phi^{\bar{\eta},0}$ using the notation in (\ref{19}). Let
$\phi_t$ and $\bar{\phi}_t$ be two solutions of the SDE (\ref{sde}) with common Brownian motions having initial data $\phi_0$ and $\bar{\phi}_0$. Let $\eta_{t}$ and
$\bar{\eta}_{t}$ be defined by $\eta_{t}(b):=\nabla\phi(b)$ and $\bar{\eta}_{t}(b):=\nabla\bar{\phi}(b)$, for all $b\in (\zd)^*$.
Since $\mu[\xi],\bar{\mu}[\xi]$ are stationary for the SDE (\ref{sde1}), we conclude that $\eta_{t}$ and
$\bar{\eta}_{t}$ are distributed by $\mu[\xi]$ and $\bar{\mu}[\xi]$ respectively, for all $t\ge 0$. 

~

\noindent We will prove

\begin{lem}
\label{annealeduniq}
For all $u\in\RR^d$, we have
\begin{equation}
\label{0}
{{\lim}}_{T\rightarrow\infty}\int\frac{1}{T}\int_0^T\sum_{b\in(\zd)^*} e^{-2r|x_b|} {\mathbb E}_{{\mathbb Q}[\xi]}\left[\left(\eta_{t}(b)-{\bar{\eta}}_{t}(b)\right)^2\right]\rmd t\P(\rmd\xi)=0.
\end{equation}
\end{lem}
By means of Proposition \ref{kom}, we will then perform an average over the integrating quantity above and find a deterministic sequence $(m_r)_{r\in\N}$, along which this average converges for $\P$-a.e. $\xi$. More precisely, we will show
\begin{lem}
\label{2.1}
There exists a deterministic sequence $(m_r)_{r\in\N}$ in $\N$ such that for $\P$-almost every $\xi$ 
\begin{equation}
\label{000a}
{{\lim}}_{k\rightarrow\infty}\frac{1}{k}\bigg(\sum_{i=1}^k\frac{1}{{m_i}}\int_0^{{m_i}}\sum_{b\in(\zd)^*} e^{-2r|x_b|}{\mathbb E}_{{\mathbb Q}[\xi]}\left[\left(\eta_{t}(b)-{\bar{\eta}}_{t}(b)\right)^2\right]\rmd t\bigg)=0.
\end{equation}
\end{lem}
Once Lemma \ref{2.1} is proved, Theorem \ref{extremeuniq} immediately follows.
Indeed Lemma \ref{2.1} implies for $\mathbb P$-almost all $\xi$
\begin{equation}
\label{j1}
\lim_{k\rightarrow\infty}\int |\eta - \bar\eta|^2_r \hat\P_k[\xi](d\eta d\bar\eta)=0,
\end{equation}
where $\hat\P_k[\xi]$ is a shift-covariant probability measure on $\chi_r\times\chi_r$,
$r>0$, defined by
$$
\hat\P_k[\xi](d\eta d\bar\eta) := \frac{1}{k}\bigg(\sum_{i=1}^k\frac{1}{{m_i}} \int^{{m_i}}_0 {\mathbb Q}[\xi](\{\eta_t(b),\bar\eta_t(b)\}_b \in d\eta d\bar\eta) ~\rmd t\bigg).
$$
The first marginal of $\hat\P_k[\xi]$ is $\mu[\xi]$ and the second one is $\bar\mu[\xi]$. Thus (\ref{j1}) implies
that the Wasserstein distance between $\mu$ and $\bar\mu$ vanishes and hence $\mu[\xi]=\bar\mu[\xi]$ for $\mathbb P$-almost all $\xi$ (see,
e.g., [14, p.482] for the Wasserstein metric on the space $\PP(\chi_r)$). This proves
Theorem \ref{extremeuniq}.

\vspace{2mm}

\textbf{Proof of Lemma \ref{2.1}}

From Proposition \ref{kom} and Lemma \ref{annealeduniq}, it follows that there exist a deterministic sequence $(m_r)_{r\in\N}$ in $\N$ and a \textit{positive} integrable random variable $X$ such that 
$${{\lim}}_{k\rightarrow\infty}\frac{1}{k}\bigg(\sum_{i=1}^k\frac{1}{{m_i}}\int_0^{{m_i}}\sum_{b\in(\zd)^*} e^{-2r|x_b|} {\mathbb E}_{{\mathbb Q}[\xi]}\left[\left(\eta_{t}(b)-{\bar{\eta}}_{t}(b)\right)^2\right]\rmd t\bigg)=X~~\mbox{for}~~\P-\mbox{almost every}~\xi.$$
It remains to show that $X=0$ for $\P$-almost every $\xi$. We note now that for all $k\ge 1$, we have
\begin{eqnarray*}
\lefteqn{\frac{1}{k}\bigg(\sum_{i=1}^k\frac{1}{{m_i}}\int_0^{{m_i}}\sum_{b\in(\zd)^*} e^{-2r|x_b|} {\mathbb E}_{{\mathbb Q}[\xi]}\left[\left(\eta_{t}(b)-{\bar{\eta}}_{t}(b)\right)^2\right]\rmd t\bigg)}\\
&\le& \frac{1}{k}\bigg(\sum_{i=1}^k\frac{2}{{m_i}}\int_0^{{m_i}}\sum_{b\in(\zd)^*} e^{-2r|x_b|} {\mathbb E}_{{\mu}[\xi]}\left(\eta_{t}(b)\right)^2\rmd t+\sum_{i=1}^k\frac{2}{{m_i}}\int_0^{{m_i}}\sum_{b\in(\zd)^*} e^{-2r|x_b|} {\mathbb E}_{{\bar\mu}[\xi]}\left(\eta_{t}(b)\right)^2\rmd t\bigg)\\
&=&2\sum_{b\in(\zd)^*} e^{-2r|x_b|} {\mathbb E}_{{\mu}[\xi]}\left(\eta(b)\right)^2+2\sum_{b\in(\zd)^*} e^{-2r|x_b|} {\mathbb E}_{{\bar\mu}[\xi]}\left(\eta(b)\right)^2,
\end{eqnarray*}
where in the equality we used that $\mu[\xi]$ and $\bar\mu[\xi]$ are stationary for the SDE (\ref{sde1}) for all fixed $\xi$. Due to the integrability assumption satisfied by $\mu[\xi]$ and $\bar\mu[\xi]$, we can now apply the Dominated Convergence Theorem to get
\begin{eqnarray*}
\E(X)&=&\E\bigg(\lim_{k\rightarrow\infty}\frac{1}{k}\bigg(\sum_{i=1}^k\frac{1}{{m_i}}\int_0^{{m_i}}\sum_{b\in(\zd)^*} e^{-2r|x_b|} {\mathbb E}_{{\mathbb Q}[\xi]}\left[\left(\eta_{t}(b)-{\bar{\eta}}_{t}(b)\right)^2\right]\rmd t\bigg)\bigg)\\
&=&\lim_{k\rightarrow}\frac{1}{k}\sum_{i=1}^k\E\bigg(\bigg(\frac{1}{{m_i}}\int_0^{{m_i}}\sum_{b\in(\zd)^*} e^{-2r|x_b|} {\mathbb E}_{{\mathbb Q}[\xi]}\left[\left(\eta_{t}(b)-{\bar{\eta}}_{t}(b)\right)^2\right]\rmd t\bigg)\bigg).
\end{eqnarray*}
Coupled with (\ref{0}), the above gives by the Ces\`aro Means theorem that $\E(X)=0$, and therefore $X=0$ for $\P$-almost every $\xi$.

$\Cox$

\vspace{2mm}

\textbf{Proof of Lemma \ref{annealeduniq}.}

We will use in our proof the following notations for the measures averaged over the disorder
$$\mu_{av}(\rmd\eta):=\left(\int\mathbb{P}(d\xi)\mu[\xi]\right)(\rmd\eta),~~{\bar{\mu}}_{av}(\rmd\bar\eta):=\left(\int\mathbb{P}(d\xi)\bar\mu[\xi]\right)(\rmd\bar\eta)~~\mbox{and}~~{\mathbb Q}_{av}:=\int \mathbb{Q}[\xi]\P(\rmd \xi).$$
We will also use in our proof the fact that $\mu[\xi]$ is stationary for the SDE (\ref{sde1}) for each fixed $\xi$.

\vspace{2mm}

By the same reasoning as in (2.10) from Proposition 2.1 in \cite{FS}, we obtain,
with the choice $\Lambda = \Lambda_\ell:=[-\ell,\ell]^d\cap\zd,\ell>0$.
\begin{multline}
\label{j2}
\E_{{\mathbb Q}[\xi]}\left[\sum_{x\in\Lambda_\ell}\left(\tilde\phi_T(x)\right)^2\right] +
C_1 \int^T_0 \E_{{\mathbb Q}[\xi]} \left[ \sum_{b\in\Lambda^*_\ell} \left(\nabla\tilde\phi_t(b)\right)^2\right] ~\rmd t\\
\le \E_{{\mathbb Q}[\xi]} \left[ \sum_{x\in\Lambda_\ell}\left(\tilde\phi_0(x)\right)^2\right]+ 2 C_2 \int^T_0 \E_{{\mathbb Q}[\xi]} \bigg[ \sum_{b\in\partial\Lambda^*_\ell\atop x_b\in\Lambda,y_b\notin\Lambda} |\tilde\phi_t(x_b)|
|\nabla \tilde\phi_t(b)|\bigg]~\rmd t,
\end{multline}
for every $T>0$ and $\ell\in\N$. We note now that the distribution of $(\eta_t,\bar\eta_t)=(\nabla\phi_t,\nabla\bar\phi_t)$ on $\chi_r\times\chi_r$ is shift-covariant due to the independence of $\eta$ and $\bar\eta$ and to the shift-covariance of $\mu[\xi]$ and $\bar\mu[\xi]$. Since the disorder is i.i.d. (respectively stationary for model B), it follows that averaging this distribution over the disorder produces a shift-invariant measure. It follows that to prove (\ref{0}), it is sufficient to show
$${{\lim}}_{T\rightarrow\infty}\frac{1}{T}\int_0^T\sum_{\alpha=1}^d\E_{{\mathbb Q}_{av}} \left(\nabla\tilde\phi_t(e_\alpha)\right)^2 ~\rmd t=0.$$
Therefore, we can now proceed as in Step 1 from \cite{FS} and we get in ({\ref{j2})
\begin{eqnarray*}
\int_0^T\sum_{\alpha=1}^d\E_{{\mathbb Q}_{av}} \left(\nabla\tilde\phi_t(e_\alpha)\right)^2 ~\rmd t\le\frac{2d}{C_1|\Lambda_l^*|}\E_{{\mathbb Q}_{av}} \left[ \sum_{x\in\Lambda_\ell}\left(\tilde\phi_0(x)\right)^2\right]+\frac{(2C_2 c_0)^2d}{(C_1l)^2}
\int_0^T\sup_{y\in\partial\Lambda_l}\|\tilde\phi_t\|^2_{{\mathbb{Q}}_{av}}\rmd t,
\end{eqnarray*}
where $c_0:=\sup_{l\ge 1}{l|\partial\Lambda^*|/|\Lambda^*|}<\infty$.

In order to use the same reasoning for our proof as in Proposition 2.1 from \cite{FS}, we need to show that a certain ergodic theorem holds for our measures averaged over the disorder. By means of the ergodic decomposition for ${{\mu}_{av}}$ there exists a probability measure $\rho_{\mu_{av}}$ on the set of ergodic measures on $\chi$, denoted by ${\cal{M}}_e(\chi)$, such that we have
$$\mu_{av}=\int_{{\cal{M}}_e(\chi)}\gamma\rho_{\mu_{av}}(\rmd\gamma).$$
In particular, for all $\alpha\in\{1,2,\ldots, d\}$, we have
$$\mu_{av}(E_{\alpha})=\int_{{\cal{M}}_e(\chi)}\gamma(E_{\alpha})\rho_{\mu_{av}}(\rmd\gamma).$$
Since by hypothesis $\mu_{av}(E_{\alpha})=1$, it follows that for all $\rho_
{\mu_{av}}$-a.e. $\gamma\in {\cal{M}}_e(\chi)$ we have $\gamma(E_{\alpha})=1$. Due to the shift-invariance of $\gamma$ this implies
$$ \gamma(\eta(b))=\langle u,y_b-x_b \rangle ~\mbox{for all bonds}~b=(x_b,y_b)\in(\Z^d)^*.$$
To bound
$$\|\phi^{\eta,0}(x) - x\cdot u\|^2_{L^2(\mu_{av})}=\int_{{\cal{M}}_e(\chi)}\gamma\left((\phi^{\eta,0}(x) - x\cdot u)^2\right)\rho_{\mu_{av}}(\rmd\gamma),$$
we will use as in \cite{FS} a special ergodic theorem for co-cycles (see for example Theorem 4 in \cite{co1}); we apply it to each $\gamma\in {\cal{M}}_e(\chi)$ to obtain
\begin{equation}
\label{ergtheo}
\lim_{|x|\rightarrow\infty}\frac{1}{|x|}\|\phi^{\eta,0}(x)-x\cdot u\|_{L^2(\gamma)}=0.
\end{equation}
Since for all $\gamma\in {\cal{M}}_e(\chi)$ 
$$\frac{1}{|x|}\|\phi^{\eta,0}(x)-x\cdot u\|^2_{L^2(\gamma)}\le \sum_{i=1}^d 2d \gamma((\eta(e_i))^2),$$
with $\sum_{i=1}^d\int_{{\cal{M}}_e(\chi)}\gamma((\eta(e_i))^2)\rmd\gamma=\sum_{i=1}^d\mu_{av}((\eta(e_i))^2)<\infty$, we have by the Dominated Convergence Theorem that
\begin{equation}
\label{ergo1}
 \lim_{|x|\rightarrow\infty} \frac{1}{|x|^2}\|\phi^{\eta,0}(x) - x\cdot u\|^2_{L^2(\mu_{av})}\le \int_{{\cal{M}}_e(\chi)}\lim_{|x|\rightarrow\infty} \frac{1}{|x|^2}\gamma\left((\phi^{\eta,0}(x) - x\cdot u)^2\right)\rho_{\mu_{av}}(\rmd\gamma)=0,
 \end{equation}
with a similar estimate holding for $\lim_{|x|\rightarrow\infty} \frac{1}{\|x\|}\|\phi^{\eta,0}(x) - x\cdot u\|^2_{L^2({\bar\mu}_{av})}$. Fix $\epsilon>0$.  It follows from (\ref{ergo1}) that there exists $l_0=l_0(\epsilon)>0$ such that for all $|x|\ge l_0$
\begin{equation}
\label{ergod1}
 \frac{1}{|x|^2}\|\phi^{\eta,0}(x) - x\cdot u\|^2_{L^2(\mu_{av})}\le\epsilon~~~\mbox{and}~~~ \frac{1}{|x|^2}\|\phi^{\eta,0}(x) - x\cdot u\|^2_{L^2(\bar\mu_{av})}\le\epsilon.
 \end{equation}
 Given (\ref{ergod1}), the proof now follows similar arguments as in \cite{FS} and will be omitted.
$\Cox$

\subsection{Ergodicity of the unique measure with given direction-averaged tilt averaged over the disorder}

In this subsection, we will show that the unique gradient measure with direction-averaged tilt $\mu[\xi]$, respectively $\mu[\omega]$, from Theorem \ref{extremeuniq} is such that the corresponding annealed measure is ergodic. We will prove 
\begin{thm}
\label{ergodic}
Let $u\in\R^d$. 
\begin{enumerate}
\item [(a)] \textbf{(Model A)} Let $d\ge 3$. Assume that $V$ satisfies (\ref{tag22}) and that $(\xi(x))_{x\in\zd}$ have symmetric distribution. For $d=3$ we will also assume that the distribution of $\xi(0)$ satisfies (\ref{gap}). Then if $\xi\rightarrow\mu[\xi]$ is the $\P$-almost surely unique shift-covariant measure $\mu[\xi]$ from Theorem \ref{extremeuniq} (a), the corresponding annealed measure ${\mu}_{av}^u(\eta):=\E\int\mu^u[\xi](\ormd\eta)$ is ergodic.

\item[(b)] \textbf{(Model B)} Let $d\ge 1$. Assume that for $\P$-almost every $\omega$, $V^\omega_{(x,y)}$ satisfies (\ref{tag23}) uniformly in the bonds $(x,y)$. Then if $\omega \rightarrow\mu[\omega]$ is the $\P$-almost surely unique shift-covariant measure $\mu[\omega]$ from Theorem \ref{extremeuniq} (b), the corresponding annealed measure ${\mu}_{av}^u(\eta):=\E\int\mu^u[\omega](\ormd\eta)$ is ergodic.

\end{enumerate}
\end{thm}

{\bf Proof.} We will only do the proof of the theorem for (a), the proof for (b) following similarly. 

\vspace{2mm}

Let ${\cal F}_{inv}(\chi)$ the $\sigma$-algebra of shift-invariant events on $\chi$ (i.e., the sets $A$ satisfying $\tau_v(A)=A$  for all $v\in\zd$). By \cite{giorgii} we need to show that for all $A\in {\cal F}_{inv}(\chi)$, we have $\mu_{av}^u(A)=0$ or $\mu_{av}^u(A)=1$. We will show that this holds by contradiction.  

Suppose that there exists $A\in {\cal F}_{inv}(\chi)$ such that $0<\mu_{av}^u(A)<1$. Then, for $\P$-almost all $\xi$ we have  $0<\mu^u[\xi](A)<1$. We define now for all $\xi$ the \textit{distinct} measures on $\chi$
$$\mu^u_A[\xi](B):=\frac{\mu^u[\xi](B\cap A)}{\mu^u[\xi](A)}~~~\mbox{and}~~~\mu^u_{A^c}[\xi](B):=\frac{\mu^u[\xi](B\cap A^c)}{\mu^u[\xi](A^c)},~~\mbox{for all}~~B\in \TT,$$
where we denoted by
$\TT:=\sigma(\{\eta_b:b\in (Z^d)^*\})$ the smallest $\sigma$-algebra on $(\zd)^*$ generated by all the edges in $(\zd)^*$. 

It is easy to show that $\mu^u_A[\xi](E_\alpha)=1$ and $\mu^u_{A^c}[\xi](E_\alpha)=1$, for $\alpha\in \{1,2,\ldots, d\}$. More precisely, in view of $\mu^u[\xi](E_\alpha)=1$, $\alpha\in \{1,2,\ldots, d\}$, we have
$$\mu^u_A[\xi](E_\alpha)=\frac{\mu^u[\xi](E_\alpha\cap A)}{\mu^u[\xi](A)}=\frac{\mu^u[\xi](E_\alpha)+\mu^u[\xi](A)-\mu^u[\xi](E_\alpha\cup A)}{\mu^u[\xi](A)}=\frac{\mu^u[\xi](A)}{\mu^u[\xi](A)}=1,$$
with a similar argument for  $\mu^u_{A^c}[\xi](E_\alpha)$.
 Moreover, since $A$ is an invariant set and $\mu^u[\xi]$ is shift-covariant, the measures $\E\int{\mu^u_{A}}[\xi](\ormd\eta)$ and $\E\int{\mu^u}_{A^c}[\xi](\ormd\eta)$ are shift-invariant. Therefore $\mu^u_A[\xi]$ and $\mu^u_{A^c}[\xi]$ satisfy all the assumptions of Theorem \ref{extremeuniq}. It follows now by Theorem \ref{extremeuniq} that $\mu^u_A[\xi]=\mu^u_{A^c}[\xi]$ for $\P$-almost all $\xi$, which leads to a contradiction. 
$\Cox$

\vspace{1mm}

As a direct consequence of Theorems  \ref{extremeuniq} and \ref{ergodic}, we get
\begin{cor}
\label{exist}
Let $u\in\R^d$. Under the assumptions of Theorem \ref{ergodic}, there exists at least one shift-covariant gradient Gibbs measure $\xi\rightarrow\mu[\xi]$ (respectively $\omega\rightarrow\mu[\omega]$) with expected given tilt $u$ and with the corresponding annealed measure being ergodic.
\end{cor}
{\bf Proof.} 


The statement follows immediately by applying Theorems \ref{extremeuniq} and \ref{ergodic}.

$\Cox$

\subsection{Proof of Theorem \ref{uniq}}

We assume that there exist at least two shift-covariant gradient Gibbs measures $\xi\rightarrow\mu[\xi]$ and $\xi\rightarrow\bar\mu[\xi]$ (respectively $\omega\rightarrow\mu[\omega]$ and $\omega\rightarrow\bar\mu[\omega]$) with expected given tilt $u$ and with the corresponding annealed measure being ergodic. By Corollary \ref{exist}, the existence of at least one such gradient Gibbs measure is assured. Due to the ergodicity of the annealed measures, (\ref{ergod1}) above holds by Theorem 4 in \cite{co1}. The proof of uniqueness follows now the same arguments as the proof of Theorem \ref{extremeuniq} above and will be omitted.

$\Cox$

\section{Decay of covariances for the annealed gradient Gibbs measure}
\label{decay1}

We will derive in this section the \textit{annealed} decay of covariances for the gradient Gibbs measure from Proposition \ref{existgibbs}. Since for lack of simple monotonicity arguments we were unable to prove that this measure is extremal for a.s. disorder, we can't make use of this fact in our computations below. We will employ in our proof the corresponding annealed covariances for the finite-volume Gibbs measures from (\ref{hat1}) (respectively from (\ref{hat2})), Proposition \ref{dgirep}, the bounds from Proposition \ref{propg} and the Poincar\'e-type inequality from (\ref{ledoux}) (which, unlike the more general inequality from Proposition \ref{golemma} does not contain a cumbersome, difficult to control, supremum in its formula).

\vspace{1.5mm}

\textbf{Proof of Theorem \ref{decay}}

\begin{enumerate}
\item [(a)]  	\textbf{Step 1:} We will show here that
\begin{equation}
\label{limitah}
\Cov(\mu^u[\xi](F(\eta)),\mu^u[\xi](G(\eta)))=\lim_{k\rightarrow\infty}\lim_{l\rightarrow\infty}\Cov(\hat\mu^u_k[\xi](F(\eta)),\hat\mu^u_l[\xi](G(\eta))),
\end{equation}
which will then allow us to use (\ref{ledoux}) to estimate, uniformly in $k,l$, the right-hand side of (\ref{limitah}).

\vspace{1mm}

Since
$$\Cov(\mu^u[\xi](F(\eta)),\mu^u[\xi](G(\eta))=\E\left(\mu^u[\xi]\left(F(\eta)-\E(\mu^u[\xi](F(\eta)))\right)\mu^u[\xi]\left(G(\eta)-\E(\mu^u[\xi](G(\eta)))\right)  \right),$$
it is sufficient to consider the case with $\E(\mu^u[\xi](F(\eta)))=\E(\mu^u[\xi](G(\eta)))=0$.
We note now that by Taylor's expansion, we have
\begin{equation}
\label{taylor}
F(\eta)=F(0)+\sum_{b\in(\zd)^*}\eta(b)\int_0^1\partial_bF(t\eta)\rmd t,
\end{equation}
where by hypothesis, the sum above is over finitely many coordinates and $\partial_b F$ is bounded for all $b\in (\zd)^*$ in the sum.  In view of (\ref{intcond}) from Proposition \ref{existgibbs} and of (\ref{taylor}), we have for $\P$-almost all $\xi$ that $\int\mu^u[\xi](\ormd\eta)F^2(\eta)<\infty$. It is now easy to show that 
\begin{equation}
\label{useless}
 \int\mu^u[\xi](\ormd\eta)F(\eta)= \lim_{k\rightarrow\infty}\int\hat\mu_k^u[\xi](\ormd\eta)F(\eta).
 \end{equation}
We will show next that $\hat\mu^u_k[\xi](F(\eta))\hat\mu^u_l[\xi](G(\eta))$ is a uniformly integrable double-sequence. Using this and (\ref{useless}), we can then apply the Vitali Convergence Theorem and obtain (\ref{limitah}). We note first that 
\begin{equation*}
\E\left(\left(\hat\mu^u_k[\xi](F(\eta))\hat\mu^u_l[\xi](G(\eta))\right)^2\right)\le\E\left( \left(\hat\mu^u_k[\xi](F(\eta))\right)^4\right)+\E\left( \left(\hat\mu^u_l[\xi](G(\eta))\right)^4\right)
\end{equation*}
It follows from the above that it suffices now to bound $\E \left(\left(\hat\mu^u_k[\xi](F(\eta))\right)^4\right)$ and $\E\left(\left(\hat\mu^u_l[\xi](G(\eta))\right)^4\right)$ uniformly in $k,l$. 
We have 
\begin{equation}
\label{rep1}
\E \left(\left(\hat\mu^u_k[\xi](F(\eta))\right)^4\right)=\Var \left(\left(\hat\mu^u_k[\xi](F(\eta))\right)^2\right)+\E^2 \left(\left(\hat\mu^u_k[\xi](F(\eta))\right)^2\right).
\end{equation}
By using (\ref{taylor}) and the assumptions on $F$, we have for some $C(F)>0$ independent of $k$ that
$$\E\left(\left(\hat\mu^u_k[\xi](F(\eta))\right)^2\right)\le C(F)\sum_{b\in (\zd)^*}\E \left(\left(\hat\mu^u_k[\xi](\left|\eta(b)\right|)\right)^2\right)\le C(F)\sum_{b\in (\zd)^*}\E \left(\hat\mu^u_k[\xi](\eta^2(b))\right).$$
By Proposition 3.6 from \cite{CK}, there exists $K>0$ such that $\sup_{k\in\N,b\in (\zd)^*}\E \left(\hat\mu^u_k[\xi](\eta^2(b))\right)<K$ so we only need to bound the variance term on the right-hand side of (\ref{rep1}) above. By (\ref{ledoux}) for the first inequality below, by $(\sum_{i\in I} a_i)^2\le |I| \sum_{i\in I} a_i^2, I\subset\zd$, for the second inequality and by Proposition \ref{dgirep} for the third inequality, we have for all $k\in\N$ with the notation $b=(x_b,y_b)$
\begin{eqnarray}
\label{soannoyed}
\lefteqn{\Var \left(\left(\hat\mu^u_k[\xi](F(\eta))\right)^2\right)}\nonumber\\
&\le&4C(d)\sum_{ z\in\zd}\int \left(\hat\mu^u_k[\xi](F(\eta))\right)^2\left(\frac{\partial\hat\mu^u_k[\xi]( F(\eta))}{\partial\xi(z)}\right)^2\rmd\P\nonumber\\
&\le&\frac{4C(d)}{k}\sum_{i=1}^k\frac{1}{|\L_{m_i}|}\sum_{w\in \L_{m_i}}\sum_{ z\in{\L_{m_i}+w}}\int \left(\hat\mu^u_k[\xi](F(\eta))\right)^2\cov^2_{\mu^{\rho_u}_{\L_{m_i}+w}[\xi]}(\phi(z), F(\eta))\rmd\P \nonumber\\
&\le&\frac{4C(d)}{k}\sum_{b\in (\zd)^*}\sum_{i=1}^k\frac{C_1(F)}{|\L_{m_i}|}\sum_{w\in \L_{m_i}}\sum_{ z\in{\L_{m_i}+w}}\int \left(\hat\mu^u_k[\xi](F(\eta))\right)^2\mu^{\rho_u}_{\L_{m_i}+w}[\xi]\left(\left(\nabla_{(x_b,y_b)}g^{\nabla\phi}_{\L_{m_i+w}}(x_b,z)\right)^2\right)\rmd\P,\nonumber\\
\end{eqnarray}
for some $C_1(F)>0$ which depends only on $F$ and for some $C(d)>0$ which depends only on $d$ and on the distribution of the disorder $\xi(0)$. We denoted in the above $\nabla_{(x_b,y_b)}g^{\nabla\phi}_{\L_{m_i+w}}(x_b,z):=g^{\nabla\phi}_{\L_{m_i+w}}(x_b,z)-g^{\nabla\phi}_{\L_{m_i}+w}(y_b,z)$. By Proposition \ref{propg} (i) (for $d\ge 5$) and (iv) (for $d=4$), we have
\begin{equation}
\label{rough}
\sup_{b\in (\zd)^*}\sum_{ z\in{\L_{m_i+w}}}\big(\nabla_{(x_b,y_b)}g^{\nabla\phi}_{\L_{m_i+w}}(x_b,z)\big)^2<\tilde{C}(d)<\infty,
\end{equation}
for some $\tilde{C}(d)>0$ which does not depend on $k, m_i, w$ and $b$. Therefore, we have from (\ref{soannoyed}) and (\ref{rough}) that
\begin{equation*}
\label{soannoyedah1}
\sup_k\Var \left(\left(\hat\mu^u_k[\xi](F(\eta))\right)^2\right) \le 4 C(d)C_1(F)\tilde{C}(d)\sup_k\int  \left(\hat\mu^u_k[\xi](F(\eta))\right)^2\rmd\P<\infty.
\end{equation*}
Thus $\sup_{k,l}\E\left(\left(\hat\mu^u_k[\xi](F(\eta))\hat\mu^u_l[\xi](G(\eta))\right)^2\right)<\infty$ for $d\ge 4$, 
so $\hat\mu^u_k[\xi](F(\eta))\hat\mu^u_l[\xi](G(\eta))$ is a uniformly integrable double-sequence and (\ref{limitah}) follows. However, we cannot argue for $d= 3$ that (\ref{rough}) holds based on the bounds from Proposition \ref{propg} unless the \textit{unknown} value $\delta$ from (\ref{nash}) in Proposition \ref{propg} (iv) would be known to be $>1/2$. Assume $\delta\le 1/2$. In this case, the argument is more delicate and we will proceed as follows after the last line of (\ref{soannoyed}). First
\begin{eqnarray}
\label{fedup3}
\lefteqn{\frac{1}{k}\sum_{i=1}^k \frac{1}{|\L_{m_i}|}\sum_{w\in\L_{m_i}}\sum_{ z\in{\L_{m_i}+w}}\int \left(\hat\mu^u_k[\xi](F(\eta))\right)^2\mu^{\rho_u}_{\L_{m_i}+w}[\xi]\left(\left(\nabla_{(x_b,y_b)}g^{\nabla\phi}_{\L_{m_i+w}}(x_b,z)\right)^2\right)\rmd\P}\nonumber\\
&\le&\sum_{i=1}^k \frac{1}{k|\L_{m_i}|}\sum_{z\in {\L_{m_i}+w}\atop w\in\L_{m_i}} \int \left|\left(\hat\mu^u_k[\xi](F(\eta))\right)^2-\E\left(\left(\hat\mu^u_k[\xi](F(\eta))\right)^2\right)\right| \mu^{\rho_u}_{\L_{m_i}+w}[\xi]\left(\left(\nabla_{(x_b,y_b)}g^{\nabla\phi}_{\L_{m_i+w}}(x_b,z)\right)^2\right)\rmd\P\nonumber\\
&&+\frac{1}{k}\sum_{i=1}^k \frac{1}{|\L_{m_i}|}\sum_{w\in\L_{m_i}}\sum_{ z\in{\L}_{m_i}+w}\E\left(\left(\hat\mu^u_k[\xi](F(\eta))\right)^2\right)\int \mu^{\rho_u}_{\L_{m_i}+w}[\xi]\left(\left(\nabla_{(x_b,y_b)}g^{\nabla\phi}_{\L_{m_i+w}}(x_b,z)\right)^2\right)\rmd\P.
\end{eqnarray}
The last term in the above can be bound uniformly in $k$ by similar arguments as the $d=3$ case from Theorem \ref{directilt}, and by using $\sup_k\E\left(\left(\hat\mu^u_k[\xi](F(\eta))\right)^2\right)<K$. 

It remains to bound the first term on the right-hand side in (\ref{fedup3}). By using $ab<\lambda a^2+{\lambda}^{-1}b^2,a,b\in\R,\lambda>0, g^{\nabla\phi}_{{\L}_{m_i}+w}(x,z)= g^{\tau_{-z}(\nabla\phi)}_{{\L}_{m_i}+w-z}(x-z,0)$ and the fact that 
$${\L}_{m_i}+w\subset\L_{2}\cup\cup_{j=1}^{1+\left[\log({3m_i})\right]} \left({\L}_{2^{j+1}}\setminus {\L}_{2^{j}}\right),~~\forall~m_i\in\N,w\in\L_{m_i},$$ 
we have for all $0<\alpha<1$ and for $\bar{C}>0$ to be chosen later
\begin{eqnarray}
\label{dogtired1}
\lefteqn{\sum_{i=1}^k \frac{1}{k|\L_{m_i}|}\sum_{z\in{\L}_{m_i}+w\atop w\in\L_{m_i}} \int \left|\left(\hat\mu^u_k[\xi](F(\eta))\right)^2-\E\left(\left(\hat\mu^u_k[\xi](F(\eta))\right)^2\right)\right| 
\mu^{\rho_u}_{\L_{m_i}+w}[\xi]\left(\left(\nabla_{(x_b,y_b)}g^{\nabla\phi}_{\L_{m_i+w}}(x_b,z)\right)^2\right)\rmd\P}\nonumber\\
&\le&\sum_{i=1}^k \frac{1}{k|\L_{m_i}|}\sum_{w\in\L_{m_i}}\sum_{j=0}^{1+[\log({3m_i})]} \sum_{z\in\Lambda_{2^{j+1}\setminus\Lambda_{2^j}}}\bigg(\bar{C}2^{-(j+1)(3+\alpha)}\Var \left(\left(\hat\mu^u_k[\xi](F(\eta))\right)^2\right)\nonumber\\
&&\,\,\,\,\,\,\,\,\,\,\,\,\,\,\,\,\,\,\,\,\,\,\,\,\,\,\,\,\,\,\,\,\,\,\,\,\,\,\,\,\,\,\,\,\,\,\,\,\,\,\,\,\,\,\,\,\,\,\,\,\,\,\,\,\,\,\,\,\,\,\,\,\,\,\,\,\,\,\,\,\,\,\,\,\,\,\,\,+2^{(j+1)(3+\alpha)}{\bar{C}}^{-1}\E_{\mu^{\rho_u}_{\L_{m_i}+w-z}[\xi]}\left(\left(\nabla_{(x_b,y_b)}g^{\nabla\phi}_{\L_{m_i+w-z}}(x_b-z,0)\right)^4\right)\bigg),\nonumber\\
\end{eqnarray} 
where by abuse of notation we have written $\Lambda_2\setminus\Lambda_1$ for the set $\Lambda_2$. We will next estimate separately each of the two terms on the right-hand side in (\ref{dogtired1}) above. The first term can be easily bound by 
\begin{equation}
\label{dogtired2}
\sum_{j=0}^\infty\frac{\bar{C}\Var \left(\left(\hat\mu^u_k[\xi](F(\eta))\right)^2\right)}{(2^{\alpha})^j}\le 2^{\alpha}\bar{C}/(2^\alpha-1)\Var \left(\left(\hat\mu^u_k[\xi](F(\eta))\right)^2\right).
\end{equation}
To bound the second term, we have by means of Lemma 2.9 from \cite{GO}
\begin{eqnarray}
\label{dogtired3}
\lefteqn{\frac{1}{|\Lambda_{m_i}|}\sum_{j=0}^{1+[\log({3m_i})]}2^{(j+1)(3+\alpha)}{\bar{C}}^{-1} \sum_{w\in\L_{m_i}}\sum_{z\in\Lambda_{2^{j+1}\setminus\Lambda_{2^j}}}\E_{\mu^{\rho_u}_{\L_{m_i}+w-z}[\xi]}\left(\left(\nabla_{(x_b,y_b)}g^{\nabla\phi}_{\L_{m_i+w-z}}(x_b-z,0)\right)^4\right)}\nonumber\\
&\le&\sum_{j=0}^{1+[\log({3m_i})]}\frac{2^{(j+1)(3+\alpha)}{\bar{C}}^{-1}}{|\L_{m_i}|} \sum_{v\in\L_{2m_i}}\sum_{w\in\L_{m_i}, z\in\Lambda_{2^{j+1}\setminus\Lambda_{2^j}}\atop w-z=v}\E_{\mu^{\rho_u}_{\L_{m_i}+w-z}[\xi]}\left(\left(\nabla_{(x_b,y_b)}g^{\nabla\phi}_{\L_{m_i+w-z}}(x_b-z,0)\right)^4\right)\nonumber\\
&\le&\frac{1}{|\L_{m_i}|}\sum_{v\in\L_{2m_i}}\sum_{j=0}^{1+[\log({3m_i})]}2^{(j+1)(3+\alpha)}{\bar{C}}^{-1} 2^{-5j}\le \bar{\bar{C}},
\end{eqnarray}
for some $\bar{\bar{C}}$ independent of $m_i$ and $k$. Choosing now $\bar{C}$ with $2^\alpha\bar{C}/(2^\alpha-1)<1$, we get from combining (\ref{soannoyed}), (\ref{fedup3}), (\ref{dogtired1}), (\ref{dogtired2}) and (\ref{dogtired3}) that $\sup_k\Var \left(\left(\hat\mu^u_k[\xi](F(\eta))\right)^2\right)<\infty$ and (\ref{limitah}) follows.

\vspace{1mm}

\textbf{Step 2:} We will bound here the term on the right-hand side of (\ref{limitah}), uniformly in $k,l\in\N$, by means of (\ref{ledoux}), Proposition \ref{dgirep} and Proposition \ref{propg}.

\vspace{1mm}

First, by means of (\ref{ledoux}) we have for all $k,l\in\N$ for some $C_5(d)>0$ depending only on $d$ and on the distribution of $\xi(0)$
\begin{eqnarray}
\label{covbound}
\lefteqn{\left| \Cov(\hat\mu^u_k[\xi](F(\eta)),\hat\mu^u_l[\xi](G(\eta)))\right|}\nonumber\\
&\le&C_5(d)\sum_{ z\in\zd}\left(\int \left(\frac{\partial\hat\mu^u_k[\xi](F(\eta))}{\partial\xi(z)}\right)^2\rmd\P\right)^{1/2}\left(\int \left(\frac{\partial\hat\mu^u_l[\xi](G(\eta))}{\partial\xi(z)}\right)^2\rmd\P\right)^{1/2}\nonumber\\
&\le& C_5(d)\sum_{z\in\L_{k,l}}\E^{1/2}\bigg[\bigg(\sum_{b\in (\zd)^*\atop b=(x_b,y_b)}\sum_{i=1}^k\frac{||\partial_b F||_\infty}{k|\L_{m_i}|}\sum_{w\in \L_{m_i}}\mu^{\rho_u}_{\L_{m_i}+w}[\xi]\left( g^{\nabla\phi}_{{\L}_{m_i}+w}(z, x_b)-g^{\nabla\phi}_{{\L}_{m_i}+w}(z,y_b)\right)\bigg)^2\bigg]\nonumber\\
&&\,\,\,\,\,\,\,\,\,\,\,\,\,\,\,\,\,\,\,\,\,\,\,\,\,\,\,\,\,\E^{1/2}\bigg[\bigg(\sum_{b'\in (\zd)^*\atop b'=(x_{b'}, y_{b'})}\sum_{j=1}^l\frac{||\partial_{b'} G||_\infty}{l|\L_{m_j}|}\sum_{v\in \L_{m_j}}\mu^{\rho_u}_{\L_{m_j}+v}[\xi]\left( g^{\nabla\phi}_{{\L}_{m_j}+v}(z, x_{b'})-g^{\nabla\phi}_{{\L}_{m_j}+v}(z, y_{b'})\right)\bigg)^2\bigg]\nonumber\\
&\le&C_5(d)\sum_{z\in\L_{k,l}}\E^{1/2}\bigg(\sum_{b\in (\zd)^*\atop b=(x_b,y_b)}\sum_{i=1}^k\frac{||\partial_b F||^2_\infty}{k|\L_{m_i}|}\sum_{w\in \L_{m_i}}\mu^{\rho_u}_{\L_{m_i}+w}[\xi]\left(g^{\nabla\phi}_{{\L}_{m_i}+w}(z, x_b)-g^{\nabla\phi}_{{\L}_{m_i}+w}(z,y_b)\right)^2\bigg)\nonumber\\
&&\,\,\,\,\,\,\,\,\,\,\,\,\,\,\,\,\,\,\,\,\,\E^{1/2}\bigg(\sum_{b'\in (\zd)^*\atop b'=(x_{b'},y_{b'})}\sum_{j=1}^l\frac{||\partial_{b'} G||^2_\infty}{l|\L_{m_j}|}\sum_{v\in \L_{m_j}}\mu^{\rho_u}_{\L_{m_j}+v}[\xi]\left(g^{\nabla\phi}_{{\L}_{m_j}+v}(z, x_{b'})-g^{\nabla\phi}_{{\L}_{m_j}+v}(z, y_{b'})\right)^2\bigg),\nonumber\\
\end{eqnarray}
where  $\L_{k,l}:=\L_{2m_{\max(k,l)}}$, the first inequality above follows by Proposition \ref{dgirep}, and for the second one we used $(\sum_{i\in I} a_i)^2\le |I| \sum_{i\in I} a_i^2, I\subset\zd$. We recall here that the sums over $b,b'\in (\zd)^*$ are finite. 
To further bound (\ref{covbound}) and obtain the \textit{optimal} covariance estimates from Theorem \ref{decay}, we need to work with the infinite-volume gradient Gibbs measure $\mu^u[\xi]$ and with the infinite-volume Green's function $g$, rather than with  the corresponding finite-volume gradient Gibbs measures and finite-volume Green's functions from (\ref{covbound}). For this purpose, we would like to use the weak convergence of $\hat\mu^u_k[\xi]$ to $\mu^u[\xi]$ and the estimates in (\ref{ahr}), so we first need to control the sums in (\ref{covbound}) above for $k,l\rightarrow\infty$. To achieve this, we will first use
\begin{multline}
\label{covcovlala}
\sum_{z\in\L_{k,l}}\E^{1/2}\bigg(\sum_{b\in (\zd)^*\atop b=(x_b,y_b)}\sum_{i=1}^k\frac{||\partial_b F||^2_\infty}{k|\L_{m_i}|}\sum_{w\in \L_{m_i}}\mu^{\rho_u}_{\L_{m_i}+w}[\xi]\left(g^{\nabla\phi}_{{\L}_{m_i}+w}(z, x_b)-g^{\nabla\phi}_{{\L}_{m_i}+w}(z,y_b)\right)^2\bigg)\\
\,\,\,\,\,\,\,\,\,\,\,\,\,\,\,\,\,\,\,\,\,\E^{1/2}\bigg(\sum_{b'\in (\zd)^*\atop b'=(x_{b'},y_{b'})}\sum_{j=1}^l\frac{||\partial_{b'} G||^2_\infty}{l|\L_{m_j}|}\sum_{v\in \L_{m_j}}\mu^{\rho_u}_{\L_{m_j}+v}[\xi]\left(g^{\nabla\phi}_{{\L}_{m_j}+v}(z, x_{b'})-g^{\nabla\phi}_{{\L}_{m_j}+v}(z, y_{b'})\right)^2\bigg),\\
\le\sum_{z\in\L_{k,l}}\bigg[\E\bigg(\sum_{b\in (\zd)^*\atop b=(x_b,y_b)}\sum_{i=1}^k\frac{||\partial_b F||^2_\infty}{k|\L_{m_i}|}\sum_{w\in \L_{m_i}}\mu^{\rho_u}_{\L_{m_i}+w-z}[\xi]\left( \left(g^{\nabla\phi}_{{\L}_{m_i}+w-z}(0, x_b-z)-g^{\nabla\phi}_{{\L}_{m_i}+w-z}(0,y_b-z)\right)^2\right)\bigg)\\
\,\,\,\,\,+\E\bigg(\sum_{b'\in (\zd)^*\atop b'=(x_{b'}, y_{b'})}\sum_{j=1}^l\frac{||\partial_{b'} G||^2_\infty}{l|\L_{m_j}|}\sum_{v\in \L_{m_j}}\mu^{\rho_u}_{\L_{m_j}+v-z}[\xi]\left( \left(g^{\nabla\phi}_{{\L}_{m_j}+v-z}(0, x_{b'}-z)-g^{\nabla\phi}_{{\L}_{m_j}+v-z}(0, y_{b'}-z)\right)^2\right)\bigg)\bigg],
\end{multline}
where for the inequality above, we used $ab<a^2+b^2,a,b\in\R$, the same change of variables as in (\ref{gogo2}) and the fact that $(\xi(x))_{x\in\zd}$ are i.i.d.. We note now that for every fixed $b=(x_b,y_b)\in (\zd)^*$ and $1\le i\le k$, we have for $|z-x_b|>R$, where $R>0$ is arbitrarily fixed 
\begin{multline}
\label{1111}
\sum_{z\in\L_{k,l}, |z-x_b|>R}\E\bigg(\frac{1}{|\L_{m_i}|}\sum_{w\in \L_{m_i}}\mu^{\rho_u}_{\L_{m_i}+w-z}[\xi]\left( \left(g^{\nabla\phi}_{{\L}_{m_i}+w-z}(0, x_b-z)-g^{\nabla\phi}_{{\L}_{m_i}+w-z}(0,y_b-z)\right)^2\right)\bigg)\\
\le \sum_{v\in\L_{2m_i}}\E\bigg(\frac{1}{|\L_{m_i}|}\sum_{w, z\in \L_{2m_i}, |z-x_b|>R\atop w-z=v}\mu^{\rho_u}_{\L_{m_i}+w-z}[\xi]\left( \left(g^{\nabla\phi}_{{\L}_{m_i}+w-z}(0, x_b-z)-g^{\nabla\phi}_{{\L}_{m_i}+w-z}(0,y_b-z)\right)^2\right)\bigg)\\
\le \frac{1}{|\L_{m_i}|}\sum_{v\in\L_{2m_i}}\E\bigg(\sum_{k=0}^{\log\left(\frac{dm_i}{R_0}\right)}\sum_{2^k R\le |z-x_b|\le 2^{k+1} R}\mu^{\rho_u}_{\L_{m_i}+v}[\xi]\left( \left(g^{\nabla\phi}_{{\L}_{m_i}+v}(0, x_b-z)-g^{\nabla\phi}_{{\L}_{m_i}+v}(0,y_b-z)\right)^2\right)\bigg)\\
\le\frac{C'(d)}{R^{d-2}},\,\,\,\,\,\,\,\,\,\,\,\,\,\,\,\,\,\,\,\,\,\,\,\,\,\,\,\,\,\,\,\,\,\,\,\,\,\,\,\,\,\,\,\,\,\,\,\,\,\,\,\,\,\,\,\,\,\,\,\,\,\,\,\,\,\,\,\,\,\,\,\,\,\,\,\,\,\,\,\,\,\,\,\,\,\,\,\,\,\,\,\,\,\,\,\,\,\,\,\,\,\,\,\,\,\,\,\,\,\,\,\,\,\,\,\,\,\,\,\,\,\,\,\,\,\,\,\,\,\,\,\,\,\,\,\,\,\,\,\,\,\,\,\,\,\,\,\,\,\,\,\,\,\,\,\,\,\,\,\,\,\,\,\,\,\,\,\,\,\,\,\,\,\,\,\,\,\,\,\,\,\,\,\,\,\,\,\,\,\,\,\,\,\,\,\,\,\,\,\,\,\,\,\,\,\,\,\,\,\,\,\,\,\,\,\,\,
\end{multline}
for some $C'(d)>0$, which depends only on $d, C_1$ and $C_2$, and where for the last inequality in the above we used (\ref{gogo30a}) from Proposition \ref{propg}, with a similar inequality holding for the term on the last line of (\ref{covcovlala}). Fix $R>0$. It follows from (\ref{covbound}), (\ref{covcovlala}), (\ref{1111}) and the fact that we sum over a finite number of $b,b'\in (\zd)^*$ that 
\begin{multline}
\label{covbound11}
\left| \Cov(\hat\mu^u_k[\xi](F(\eta)),\hat\mu^u_l[\xi](G(\eta))\right|\\
\le C_5(d)\sum_{z:\max_b |z-x_b|<R\atop\max_{b'}|z-x_{b'}|<R}\E^{1/2}\bigg(\sum_{b\in (\zd)^*\atop b=(x_b,y_b)}\sum_{i=1}^k\frac{\|\partial_b F\|^2_\infty}{k|\L_{m_i}|}\sum_{w\in \L_{m_i}}\mu^{\rho_u}_{\L_{m_i}+w}[\xi]\left(g^{\nabla\phi}_{{\L}_{m_i}+w}(z, x_b)-g^{\nabla\phi}_{{\L}_{m_i}+w}(z,y_b)\right)^2\bigg)\\
\,\,\,\,\,\,\,\,\,\,\,\,\,\,\,\,\,\,\,\,\,\E^{1/2}\bigg(\sum_{b'\in (\zd)^*\atop b'=(x_{b'},y_{b'})}\sum_{j=1}^l\frac{||\partial_{b'} G||^2_\infty}{l|\L_{m_j}|}\sum_{v\in \L_{m_j}}\mu^{\rho_u}_{\L_{m_j}+v}[\xi]\left(g^{\nabla\phi}_{{\L}_{m_j}+v}(z, x_{b'})-g^{\nabla\phi}_{{\L}_{m_j}+v}(z, y_{b'})\right)^2\bigg)+\frac{C'(d)}{R^{d-2}}\\
\le C_5(d)\sum_{z:\max_{b} |z-x_b|<R\atop \max_{b'}|z-x_{b'}|<R}\E^{1/2}\bigg(\sum_{b\in (\zd)^*\atop b=(x_b,y_b)}\|\partial_b F\|^2_\infty\hat\mu^u_k[\xi]\left(g^{\nabla\phi}(z, x_b)-g^{\nabla\phi}(z,y_b)\right)^2\bigg)\\
\E^{1/2}\bigg(\sum_{b'\in (\zd)^*\atop b'=(x_{b'},y_{b'})}\|\partial_{b'} G\|^2_\infty\hat\mu^u_l[\xi]\left(g^{\nabla\phi}(z, x_{b'})-g^{\nabla\phi}(z,y_{b'})\right)^2\bigg)+\frac{C'(d)}{R^{d-2}},
\end{multline}
for some $C''(d)>0$ which depends only on $d, C_1$ and $C_2$. We used for the second inequality above the following reasoning: $g^{\nabla\phi}$ depends on $\nabla\phi$ only through $C_1\le a^{\nabla\phi}\le C_2$, from which $g^{\nabla\phi}_{{\L}_{m_i}+w}(z, x_b)-g^{\nabla\phi}_{{\L}_{m_i}+w}(z,y_b)$ converges to $g^{\nabla\phi}(z, x_b)-g^{\nabla\phi}(z,y_b)$ uniformly in $\nabla\phi$. Since the sums above are after a finite number of $z, b,b'$, we can now take limits for the finite-volume Green's functions under the expectations in the first inequality above. (To prove the uniform convergence, we apply Dini's theorem for uniform convergence: $[C_1,C_2]^{\chi}$ is compact in the product topology by Tychonoff's theorem, $\L_N\rightarrow g^{\cdot}_{\L_N}(z,x_b)$ is a non-decreasing sequence of continuous functions and the limit $g^{\cdot}(z,x_b)$ is also continuous; moreover, for all $w\in\L_{m_i}$ we have $g^{\nabla\phi}_{[0,\pm m_i]\times\ldots\times [0,\pm m_i]}(z, x_b)\le g^{\nabla\phi}_{{\L}_{m_i}+w}(z, x_b)\le g^{\nabla\phi}_{{\L}_{2m_i}}(z, x_b)$, with the sign of each $m_i$ in the lower bound interval product $[0,\pm m_i]\times\ldots\times [0,\pm m_i]$ depending on the sign of the corresponding coordinate in $w$). From (\ref{covbound}) and (\ref{covbound11}), we get
\begin{eqnarray}
\label{game}
\lefteqn{\lim_{k\rightarrow\infty}\lim_{l\rightarrow\infty}\Cov(\hat\mu^u_k[\xi](F(\eta)),\hat\mu^u_l[\xi](G(\eta))}\nonumber\\
&\le&C_5(d)\sum_{b,b'\in (\zd)^*, b=(x_b,y_b)\atop b'=(x_{b'},y_{b'})}\|\partial_b F\|_\infty\|\partial_{b'} G\|_\infty\sum_{z\in\zd}\bigg\{\E^{1/2}\bigg(\mu^u[\xi]\left(g^{\nabla\phi}(z, x_b)-g^{\nabla\phi}(z,y_b)\right)^2\bigg)\nonumber\\
&&\,\,\,\,\,\,\,\,\,\,\,\,\,\,\,\,\,\times\E^{1/2}\bigg(\mu^u[\xi]\left(g^{\nabla\phi}(z, x_{b'})-g^{\nabla\phi}(z,y_{b'})\right)^2\bigg)\bigg\},
\end{eqnarray}
where for the above we used  in the last inequality in (\ref{covbound11}) the weak convergence of $\hat\mu^u_k[\xi]$ and of $\hat\mu^u_l[\xi]$ to $\mu^u[\xi]$ (which hold in (\ref{covbound11}) since we are only summing after $z$ such that $|z-x_b|<R, |z-x_{b'}|<R$, and we are summing after a finite number of $b,b'\in (\zd)^*$) and then we took $R\rightarrow 0$. 

Given that $\E\mu^u[\xi]$ is a shift-invariant measure, we obtain now in (\ref{game}) by Proposition \ref{propg} (v)
\begin{eqnarray*}
\Cov(\mu^u[\xi](F(\eta)),\mu^u[\xi](G(\eta))\le C_5(d)\sum_{b,b'\in (\zd)^*, b=(x_b,y_b)\atop b'=(x_{b'},y_{b'})}\|\partial_b F\|_\infty\|\partial_{b'} G\|_\infty\sum_{z\in\zd}\frac{1}{]|z-x_b|[^{d-1}]|z-x_{b'}|[^{d-1}}.
\end{eqnarray*}
The statement of the theorem follows now from (\ref{ineq1}) in Proposition \ref{ineq} below.

\item [(b)] We first need to show that
\begin{equation}
\label{limitagain}
\Cov(\mu^u[\omega](F(\eta)),\mu^u[\omega](G(\eta)))=\lim_{k\rightarrow\infty}\lim_{l\rightarrow\infty}\Cov(\hat\mu^u_k[\omega](F(\eta)),\hat\mu^u_l[\omega](G(\eta)))
\end{equation}
holds. We note first that by using (\ref{taylor}) and the assumptions on $F, G$, we have for some $C(F,G)>0$ independent of $k,l$
\begin{multline*}
\E\left(\left(\hat\mu^u_k[\omega](F(\eta))\hat\mu^u_l[\omega](G(\eta))\right)^2\right)\le\E\left( \left(\hat\mu^u_k[\omega](F(\eta))\right)^4\right)+\E\left( \left(\hat\mu^u_l[\omega](G(\eta))\right)^4\right)\\
\le C(F,G)\sum_b\left\{\E \left(\left(\hat\mu^u_k[\omega](\left|\eta(b)\right|)\right)^4\right)+\E \left(\left(\hat\mu^u_l[\omega](\left|\eta(b)\right|)\right)^4\right)\right\}+F^4(0)+G^4(0).\end{multline*}
It follows from the above that it suffices now to bound $\E \left(\left(\hat\mu^u_k[\omega](\left|\eta(b)\right|)\right)^4\right)$ and $\E\left(\left(\hat\mu^u_l[\omega](\left|\eta(b)\right|)\right)^4\right)$ uniformly in $k,l$. This will prove the uniform integrability of the double-sequence $\hat\mu^u_k[\omega](F(\eta))\hat\mu^u_l[\omega](G(\eta))$, and consequently the convergence in (\ref{limitagain}). However, the situation is simpler in this case than in (a) since, as explained in Theorem \ref{directilt} (b), we have $\mu^{\rho_u}_{\L_{m_i}+w}[\omega](\phi(x)-\phi(x+e_\alpha))=u_\alpha$ for all $\alpha\in\{1,2,\ldots, d\}$, $1\le i\le k$, and for all $w\in\L_{m_i}$. Therefore, by the Brascamp-Lieb inequality (\ref{111}) applied to the convex function $L(s)=|s|$ and to each $\mu^{\rho_u}_{\L_{m_i}+w}[\omega]$, we have for all $k\ge 1$
\begin{eqnarray*}
\hat\mu^u_k[\omega](\left|\eta(b)\right|)&=&\frac{1}{k}\sum_{i=1}^k\frac{1}{|\L_{m_i}|}\sum_{w\in \L_{m_i}}\mu^{\rho_u}_{\L_{m_i}+w}[\omega]\left(\left|\eta(b)\right|\right)\\
&\le&\frac{1}{k}\sum_{i=1}^k\frac{1}{|\L_{m_i}|}\sum_{w\in \L_{m_i}}\left\{\mu^{\rho_u}_{\L_{m_i}+w}[\omega]\left(\left|\eta(b)-\mu^{\rho_u}_{\L_{m_i}+w}[\omega](\eta(b))\right|\right)+\left|\mu^{\rho_u}_{\L_{m_i}+w}[\omega](\eta(b))\right|\right\} \\
&\le& C'(d)<\infty,
\end{eqnarray*}
for some $C'(d)>0$ which depends only on $d, C_1, C_2$ and $u$. Hence (\ref{limitagain}) is proved.

We proceed next as in Step 2 from (a) above to bound the right-hand side of (\ref{limitagain}), uniformly in $k,l$. For simplicity of calculations, we assume $f_{2,b}\equiv 0$ for all $b\in (\zd)^*$. Firstly, by (\ref{ledoux}) we have
 \begin{eqnarray}
\label{covlalabound}
\lefteqn{\left| \Cov(\hat\mu^u_k[\omega](F(\eta)),\hat\mu^u_l[\omega](G(\eta)))\right|}\nonumber\\
&\le& C(d)\sum_{ b\in (\zd)^*}\left(\int \bigg(\frac{\partial \hat\mu_k^{u}[\omega](F(\eta)}{\partial\omega(b)}\bigg)^2\rmd\P\right)^{1/2}\left(\int \bigg(\frac{\partial \hat\mu_l^{u}[\omega](G(\eta)}{\partial\omega(b)}\bigg)^2\rmd\P\right)^{1/2},
\end{eqnarray}
for some $C(d)$ which depends only on $d$ and on the distribution of $V_{(x,y)}^\omega(0)$. In order to estimate the above further, we need to estimate $\bigg(\frac{\partial \hat\mu_k^{u}[\omega](F(\eta)}{\partial\omega(b)}\bigg)^2$ for all $b\in (\zd)^*$. By Proposition \ref{dgirep} for the first inequality below, Cauchy-Schwarz inequality for the second inequality, and for the third inequality by use of the Brascamp-Lieb inequality and of the fact that $\mu^{\rho_u}_{\L_{m_i}+w}[\omega](\phi(x)-\phi(x+e_\alpha))=u_\alpha$ for all $\alpha\in\{1,2,\ldots, d\}$, we have for all $b=(x_b,y_b)$ and for all $k\in\N$ 
\begin{eqnarray*}
\lefteqn{\bigg(\frac{\partial \hat\mu_k^{u}[\omega](F(\eta)}{\partial\omega(b)}\bigg)^2}\\
&=&\bigg(\frac{1}{k}\sum_{i=1}^k\frac{1}{|\L_{m_i}|}\sum_{w\in\L_{m_i}}\cov_{\mu_{\L_{m_i}+w}^{\rho_u}[\omega]}\bigg(\frac{\partial V_{(x_b,y_b)}^\omega(\phi(x_b)-\phi(y_b))}{\partial\omega(b)},F(\eta)\bigg)\bigg)^2=\\
&\le&\sum_{b'=(x_{b'},y_{b'})} ||\partial_{b'} F||^2_\infty\sum_{i=1}^k\frac{1}{k|\L_{m_i}|}\sum_{w\in {\L}_{m_i}}\\
&&\left(\mu^{\rho_u}_{\L_{m_i}+w}[\omega]\left(f_1(\omega)\left|\eta(b)\right|\big| g^{\nabla\phi}_{{\L_{m_i}}+w}(x_{b'}, x_b)-g^{\nabla\phi}_{{\L_{m_i}}+w}(x_{b'},y_b)-g^{\nabla\phi}_{{\L_{m_i}}+w}(y_{b'}, x_b)+g^{\nabla\phi}_{{\L_{m_i}}+w}(y_{b'},y_b)\big|\right)\right)^2\\
&\le&\sum_{b'=(x_{b'},y_{b'})} ||\partial_{b'} F||^2_\infty\sum_{i=1}^k\frac{f_{1,b}^2(\omega)}{k|\L_{m_i}|}\sum_{w\in {\L}_{m_i}}\mu^{\rho_u}_{\L_{m_i}+w}[\omega]\left(\eta^2(b)\right)\\
&&\,\,\,\,\,\,\,\,\,\mu^{\rho_u}_{\L_{m_i}+w}[\omega]\left(\big( g^{\nabla\phi}_{{\L_{m_i}+w}}(x_{b'}, x_b)-g^{\nabla\phi}_{{\L_{m_i}+w}}(x_{b'},y_b)-g^{\nabla\phi}_{{\L_{m_i}}+w}(y_{b'}, x_b)+g^{\nabla\phi}_{{\L_{m_i}}+w}(y_{b'},y_b)\big)^2\right)\\
&\le& \tilde{C}(d)\sum_{b'} ||\partial_{b'} F||^2_\infty\sum_{i=1}^k\frac{f^2_{1,b}(\omega)}{k|\L_{m_i}|}\\
&&\,\,\,\,\,\,\,\,\,\,\,\,\,\,\,\,\,\,\,\,\,\,\,\,\,\sum_{w\in {\L}_{m_i}}\mu^{\rho_u}_{\L_{m_i}+w}[\omega]\left(\big( g^{\nabla\phi}_{{\L}_{m_i}+w}(x_{b'}, x_b)-g^{\nabla\phi}_{{\L}_{m_i}+w}(x_{b'},y_b)-g^{\nabla\phi}_{{\L}_{m_i}+w}(y_{b'}, x_b)+g^{\nabla\phi}_{{\L}_{m_i}+w}(y_{b'},y_b)\big)^2\right),
\end{eqnarray*}
for some $\tilde{C}(d)>$ which depends only on $C_1,C_2,d$ and $u$. We use next (\ref{gogo30}), Proposition \ref{propg} (v), a similar reasoning as in part (a) above, (\ref{covlalabound}) and the above bounds, to obtain
\begin{eqnarray*}
\Cov(\mu^u[\omega](F(\eta)),\mu^u[\omega](G(\eta))\le C''(d)\sum_{b\in (\zd)^*\atop b=(x_b,y_b)}\sum_{b'\in (\zd)^*\atop b'=(x_{b'},y_{b'})}||\partial_b F||_\infty||\partial_{b'} G||_\infty\sum_{z\in\zd}\frac{1}{]|z-x_b|[^{d}]|z-x_{b'}|[^{d}}.
\end{eqnarray*}
The assertion follows now from (\ref{ineq2}) in Proposition \ref{ineq} below.

\end{enumerate}

\section{Appendix}

We will state in the next Proposition inequalities (\ref{ineq1}) and (\ref{ineq2}), used in the proof of Theorem \ref{decay}. The proof follows the same arguments as Proposition A.1 from \cite{MO} and will be omitted.
\begin{prop}
\label{ineq}
Let $x,z\in\zd$. 
\begin{itemize}
\item [(a)] For all $d\ge 3$, we have for some $C(d)>0$ which depends only on $d$
\begin{equation}
\label{ineq1}
\sum_{y\in\zd}\frac{1}{]|x-y|[^{d-1}]|z-y|[^{d-1}}\le \frac{C(d)}{]|x-z|[^{d-2}}.
\end{equation}
\item [(b)] For all $d\ge 1$ we have for some $C'(d)>0$ which depends only on $d$
\begin{equation}
\label{ineq2}
\sum_{y\in\zd}\frac{1}{]|x-y|[^{d}]|z-y|[^{d}}\le \frac{C'(d)}{]|x-z|[^{d}}.
\end{equation}
\end{itemize}
\end{prop}



\section*{Acknowledgements}

C.C. thanks Jean-Dominique Deuschel for bringing \cite{GO} to her attention, and Tadahisa Funaki, Antoine Gloria and Tom Spencer for fruitful discussions. C.C. also thanks her mother, Aurelia Cotar, for teaching her how to face with true courage impossible odds.

\end{document}